\documentclass{article}

\usepackage{arxiv}
\usepackage{float}
\usepackage[utf8]{inputenc} 
\usepackage[T1]{fontenc}    
\usepackage{hyperref}       
\usepackage{url}            
\usepackage{booktabs}       
\usepackage{amsfonts}       
\usepackage{nicefrac}       
\usepackage{microtype}      
\usepackage{lipsum}
\usepackage{enumitem}
\usepackage{graphicx}
\usepackage{amsmath}
\usepackage[pagewise]{lineno}
\graphicspath{ {./images/} }

\title{Bounds related to the edge-list chromatic and total chromatic numbers of a simple graph}

\author{
  M. Henderson \\
  \texttt{matthew.james.henderson@gmail.com} \\
   \And
A.J.W. Hilton \\
  Department of Mathematics\\
  University of Reading\\
  Whiteknights\\
  Reading, RG6 6AX, England, UK\\
  Also: Department of Mathematics\\ 
  Queen Mary University of London\\
  Mile End Road\\ 
  London E1 4NS, UK\\
\texttt{A.J.W.Hilton@reading.ac.uk}
  \And
 R. Mary Jeya Jothi\\
Department of Mathematics\\ 
Sathyabama Institute of Science and Technology\\
Chennai, Tamilnadu, India\\
\texttt{rmaryjeyajothi@gmail.com}
}

\begin{document}
\maketitle
\begin{abstract}
We show that for a simple graph $G$, $c'(G)\leq\Delta(G)+2$ where $c'(G)$ is the choice index (or edge-list chromatic number) of $G$, and $\Delta(G)$ is the maximum degree of $G$.

As a simple corollary of this result, we show that the total chromatic number $\chi_T(G)$ of a simple graph satisfies the inequality $\chi_T(G)\leq\ \Delta(G)+4$ and the total choice number $c_T(G)$ also satisfies this inequality.

We also relate these bounds to the Hall index and the Hall condition index of a simple graph, and to the total Hall number and the total Hall condition number of a simple graph.
\end{abstract}


\section{Introduction}
\setlength{\parindent}{1cm} We prove a new bound for the choice index $c'(G)$ of a simple graph in the first part of this paper. In sections 5 and 6, we show how our bound for the choice index implies new bound for the total chromatic number of a graph and for the total list chromatic number (or the total choice index) of a simple graph.  Although the bounds themselves are new, the argument to get the bounds from the choice index is well-known.\par

In Section 7 we discuss the Hall index and the Hall condition index of a simple graph, and in Section 8 we discuss the total Hall number, and the total Hall condition number of a simple graph.

Suppose that $G$ is a finite simple graph, $\zeta$ is an infinite collection of colours, and $2^\zeta$ is the collection of finite subsets of $\zeta$.  An {\it edge-list assignment to G}, or list assignment to $E(G)$, is a function :$E(G)\to2^\zeta$ . If $L$ is a list assignment to $E(G)$, a {\it proper L- colouring of $E(G)$} is a function $\phi:E(G)\to \zeta$  satisfying
\begin{enumerate}[label=(\roman*)]
\item $\phi(e)\in L(e)$,
\item if $e$, $f \in E(G)$ and $e$ and $f$ have a vertex in common, then $\phi(e)\ne \phi(f)$.
\end{enumerate}

The {\it choice index or edge chromatic number $c'(G)$} is the least number $n_0$ such that whenever $L$ is list assignment to $E(G)$ with $|L(e)|\ge n_0$ $\forall e\in E(G)$, then there exists a proper L-colouring of  $E(G)$.

In the case when $L(e) = L(f)$ for $e$, $f \in E(G)$ (so the lists are all the same), then $n_0$ is {\it the chromatic index}, or {\it edge chromatic number} of $G$ and is denoted by $\chi'(G)$.

Let $\Delta(G)$ be the maximum degree in $G$. A famous result of Vizing \cite{VGVizing1964} in 1964 states: 

\noindent {\it Theorem 1. If $G$ is a finite simple graph, then}
\begin{equation*}
\Delta(G)\le \chi'(G)\le \Delta(G)+1.
\end{equation*}

The question as to the value of $c'(G)$  seems to have erupted spontaneously in the early 1980’s and there are several authors variously associated with the {\it edge-list-colouring conjecture} (see \cite{BBollobs1985List},\cite{PErds1979},\cite{RHggkvist1992})

The strongest form of the edge list colouring conjecture is

\noindent {\it Conjecture 2. For a finite multigraph $G$, without loops,}
\begin{equation*}
c'(G)=\chi'(G).
\end{equation*}

In 1995, Galvin \cite{FGalvin1995} showed that Conjecture 2 is true for bipartite multigraphs:

\noindent {\it Theorem 3. If $G$ is a bipartite multigraph then}
\begin{equation*}
c'(G)=\chi'(G)=\Delta(G).
\end{equation*}
\noindent In Theorem 6 (below), we explain another exact result along the lines of Conjecture 2.

A slightly weaker conjecture than Conjecture 2, this time for simple graphs, is:

\noindent {\it Conjecture 4. For a finite simple graph $G$},
\begin{equation*}
\Delta(G)\le \chi'(G)\le c'(G)\le\Delta(G)+1.
\end{equation*}

\noindent In this paper we prove the following weaker variant of Conjecture 2 or Conjecture 4:

\noindent {\it Theorem 5. For finite simple graph G}
\begin{equation*}
c'(G)\le \Delta(G)+2.
\end{equation*}

It may be that our proof could be improved upon so as to prove Conjecture 4. We give some thoughts about this possibility in Section 4.

We now give some preliminary definitions needed to explain the further exact result referred to above. In \cite{AJWHilton1996} Hilton defined a 2-{\it improper edge list colouring} to be an L-colouring satisfying
\begin{enumerate}[label=(\roman*)]
\item $\phi(e)\in L(e)$$\hspace{10pt}(\forall \ e  \in E(G))$,
\item $|\{e:v \in e$ and $\phi(e)=c\}|\le2$ \hspace{1cm}$(\forall \ v \in V(G), \forall \ e \in E(G), \forall \ c \in \zeta(G))$.\\

\end{enumerate}

\noindent (Thus at any vertex there can be no more than two edges of the same colour).

Let $c'_{2}(G)$, be the least number $n_0$ such that if $|L(e)|\ge n_0$ $(\forall e \in E(G))$, then there exists a 2-improper L-colouring of $E(G)$. If $L(e)=L(f)$ for all $e$, $f \in E(G)$ then we have a {\it 2-improper edge-colouring} of $G$. In this case, $n_0$  is the 2-{\it improper chromatic index} of $G$, denoted $\chi'_{2}(G)$.

In \cite{AJWHilton1996} Hilton deduced from Galvin’s theorem:

\noindent {\it Theorem 6. For a multigraph G,}
\begin{equation*}
c'_{2}(G)=\chi'_{2}(G)= \Bigg\lceil\frac {\Delta (G)}{2}\Bigg\rceil.
\end{equation*}

In \cite{AJWHilton1996} it is also shown that an analogous result holds when 2(in $c'_{2}(G), \chi'_{2}(G)$)  is replaced by any even integer, and that we may permit $G$ to have loops, these counting 2 to the degree of the vertex they are on.

\section{Colour interchange paths (CIP's)}

Recall that a path is a sequence of distinct vertices and edges, $p_1,e_1,p_2,e_2,p_3,...,p_{r-1},e_{r-1},p_r,$ where $e_i$ is incident with $p_i$ and $p_{i+1}$ $(1\leq i \leq r-1)$ and $p_1,...,p_r$ are distinct. We could also denote the path by $p_1,p_1p_2,p_2,p_2p_3,p_3,...,p_{r-1},p_{r-1}p_r,p_r$ or more conveniently $p_1p_2,p_2p_3,...,p_{r-1}p_r$.

The proof of Vizing’s theorem depends partly on colour interchange paths (CIP’s). Consider the case when we have a properly edge-coloured graph, and the set of colours available to colour an edge is the same for each edge. Suppose we have a vertex $v$ where one colour, say $\alpha$, is absent (i.e. not used on any edge incident with $v$), and another colour, say $\beta$ is present (i.e. occurs on an edge $vu$, say). There is necessarily a path whose edges are coloured alternatively $\beta$ and $\alpha$ starting at $v$ and finishing at some vertex $w$ where at least one of $\alpha$ and $\beta$  is absent. We may interchange the colours $\alpha$ and $\beta$ on this path, producing a different proper edge colouring of $G$.

Now consider the case where we have a properly edge coloured graph, but this time the lists $L(e)$ of colours available to colour each edge $e$ vary, so that now $L(e)$ need not equal $L(f)$ when $e\ne f$. We can still have colour interchange paths in this situation, but the colour interchange paths would normally involve more than two colours. Let $v$ be any vertex in a properly L-coloured graph $G$ and $p_1$ be any vertex adjacent to $v$. Now let $a_1$ be a colour in $L(vp_1)$ (the list associated with the edge $vp_1$) missing at $v$ but present at $p_1$. 

We shall suppose there is a finite path beginning at $p_1$, say $p_1 p_2,p_2p_3,...,p_{s-1}p_s$ respectively which is coloured  $a_1,a_2,...,a_{s-1}$ and that, in addition if $p_1p_2,p_2 p_3,...,p_{s-1} p_s$  were to be recoloured $a_2, a_3,...,a_s$,respectively then we would obtain a different proper L-colouring of $G$. Note that, implicit in this assumption, is the condition that there is no edge of the form $p_{s}z$ with $a_s \in L(p_{s}z)$ and $z \in \{p_1,…,p_s\}$. Note that, before the recolouring, $a_s$ was missing at $p_s$, and that after the recolouring, $a_{s-1}$ is missing at $p_s$. This is illustrated in Figure 1.

\begin{figure}[H] 
    \centering
    \includegraphics[width=12cm]{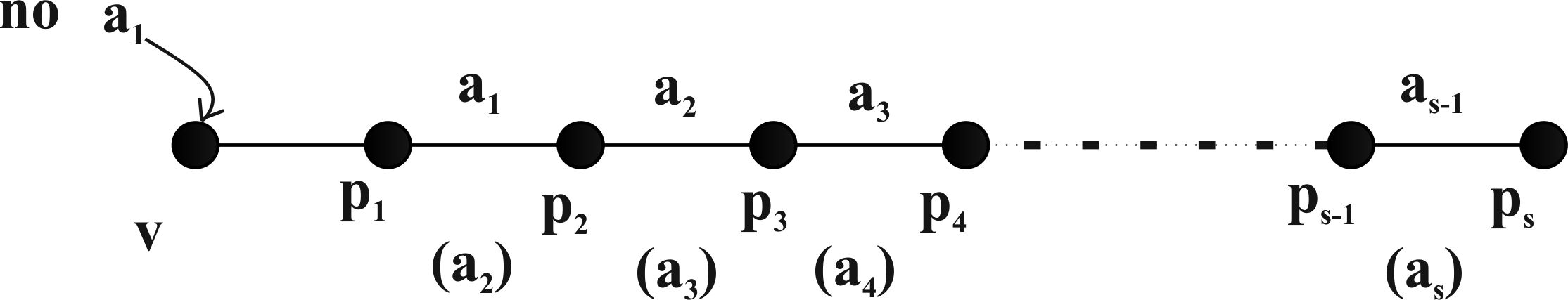}
    \caption{A Colour Interchange Path (CIP), $p_1p_2,p_2 p_3,...,p_{s-1}p_s$.}
\end{figure}

We call such a path $p_1p_2,p_2p_3,...,p_{s-1}p_s$ {\it an edge-colour interchange path starting on the edge} $p_1p_2$. We shorten this to CIP. Note that after the interchange of colours, we may recolour the edge $vp_1$ with colour $a_1$.

Let us suppose that $|L(e)|\ge \Delta(G)+2$ for all $e \in E(G)$. Suppose we have an L-edge-colouring of $G-e$ for some edge $e$ of $G$, with $e=vp_1$. Let $p_1p_2,p_2 p_3,...,p_{s-1}p_s$ be a colour interchange path starting on $p_1p_2$. At the vertices $v$ and $p_1$, in $G-e$ at least three colours are absent, and at $p_2$, at least two colours are absent. The edge $p_1p_2$  which is coloured $a_1$ satisfies $|L(p_1p_2)|\ge \Delta(G)+2\ge d_G(p_2)+2$, so there is a choice of at least two (actually three) colours in $L(p_1p_2)$ missing at $p_1$, say $a_{21}$ and $a_{22}$. Suppose we choose to place $a_{21}$ on the edge $p_1p_2$. Put $a_{21}=a_2$. Then there might be an edge coloured $a_2$ incident with $p_2$, say $p_2p_3$. Then in $L(p_2p_3)$, since $|L(p_2p_3)|\ge \Delta(G)+2$, there is a choice of at least two colours, say $a_{31}$ and $a_{32}$, not used to colour any edge incident with $p_2$. We could choose either to place on $p_2p_3$. Suppose we choose $a_{31}$, and let $a_3=a_{31}$. There might be an edge, say $p_3p_4$ coloured $a_3$. We continue in this way until the process stops. It will stop when we get to an edge $p_{s-1}p_s$ coloured $a_{s-1}$ with $a_s \in L(p_{s-1}p_s)$ with the property that no edge coloured $a_s$ is incident with $p_s$. There is no {\it a priori} reason why this process should stop, but part of our proof is a demonstration that there always is a finite colour interchange path (CIP) starting on any edge $p_1p_2$.  Since there usually is a choice of at least two colours at each step to continue constructing the path with, it should not come as a surprise to find that there always is such a CIP (if we were to assume that $|L(e)|\le \Delta(G)+1$ $(\forall e \in E(G))$, then it would not seem to be quite so likely that there is a CIP).

Intimately related to our proof of Theorem 5 is the following theorem.

\noindent {\it Theorem 7. Let $G$ be a finite simple graph . Let $E(G)$ be given a list assignment $L$. Suppose that $L$ has the property that $|L(e)|\ge \Delta(G)+2$ $(\forall e \in E(G))$. Then at each vertex $v$, each edge incident with $v$ is the start of a finite } CIP.

We may expect normally that if we have constructed part of a potential CIP, say $p_1p_2,p_2p_3,...,p_{k-1}p_k$, then there will be at least two possible edges with which to continue the potential CIP, say $p_1p_2,...,p_{k-1}p_k,p_kp_{k+1}^{'}$ and $p_1p_2,...,p_{k-1}p_k,p_kp_{k+1}^{*}$, where $p_{k+1}^{'}$ $\neq$ $p_{k+1}^{*}$. However there may be situations where for some reason we wish to restrict the choice  (so for example we might wish to allow $p_1p_2,...,p_{k-1}p_k,p_kp_{k+1}^{'}$ only, and disallow  $p_1p_2,...,p_{k-1}p_k,p_kp_{k+1}^{*}$). We call such a vertex $p_k$ in a CIP a {\it restricted vertex}. Moreover the restricted vertices in any CIP will all be the neighbours (in $G$) of some vertex, say $w$. We call a CIP {\it $w$-restricted} if the restricted vertices are all the neighbours of a vertex $w$ and we call such vertices {\it $w$-restricted vertices}. We shall show that, starting on an edge $p_1p_2$, there is a $w$-restricted finite CIP for any choice of vertex $w\notin{\{p_1,p_2\}}$.  Then we show in Part 2 that $G$ has a $w$-restricted finite CIP starting on any vertex.

\noindent We sum up the main features of w-restricted CIP’s in the following lemma.

\noindent {\it Lemma 8. In a simple graph $G$, let $p_1p_2,p_2p_3,...,p_{s-1}p_s$ be a $w$-restricted} CIP {\it $P$ with colours $a_1,a_2,...,a_s$ such that if $p_1p_2,...,p_{s-1}p_s$ are coloured $a_1,...,a_{s-1}$ respectively then this is part of a proper} L{\it-colouring of $G$, and if they are coloured $a_2,...,a_s$ respectively, then this also is part of an} L{\it -colouring of $G$}.

\noindent {\it Let $p_k$ be a $w$-restricted vertex with incident edges $p_{k-1}p_k$ and $p_kp_{k+1}$ in $P$. Then}

\begin{enumerate}[label=(\roman*)]
\item $p_k$ is {\it a neighbour} of $w$,
\item $a_{k-1}\in L(p_{k-1}p_k)$,
\item {\it There is a colour} $c\in L(p_{k-1}p_k)\cap L(p_k,w),c\neq a_{k-1}$.
\end{enumerate}

Lemma 8 is illustrated in Figure 2.

\begin{figure}[H] 
    \centering
    \includegraphics[width=11.5cm]{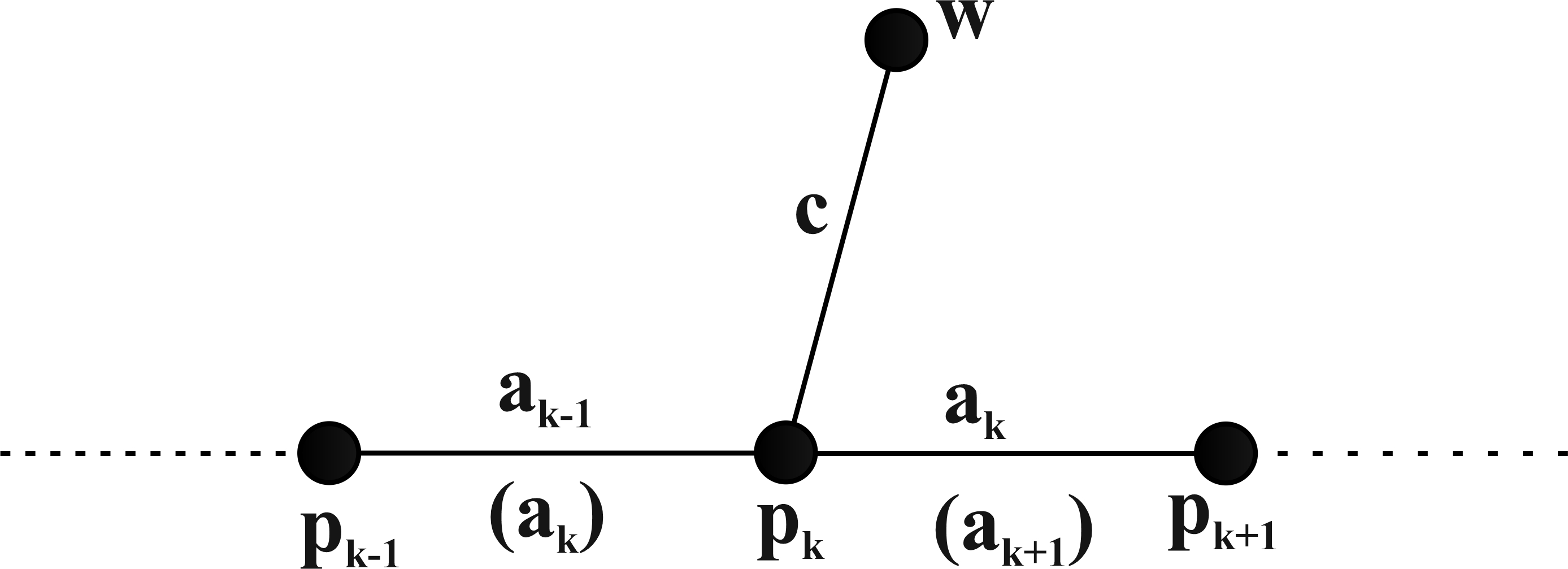}
    \caption{{\it A restricted vertex $p_k$.}}
\end{figure}

\noindent The main point here is that for some reason we may not wish to extend the path $p_1p_2,...,p_{k-1}p_k$, so that the next edge is $p_{k}w$, but prefer to extend it so that the next edge is $p_k p_{k+1}\neq p_{k}w$. 

\section{Cutting a CIP}

\noindent In the case when a colour interchange path (CIP) employs more than two colours, it is possible to ‘cut’ a CIP into two parts, both of which are CIP’s in their own right. We shall use this in Part 2 of the proof of our main theorem, Theorem 10. 

\noindent Let $G$ be a simple graph, let $P: p_1p_2,p_2p_3,...,p_{s-1}p_s$ be a $w$-restricted CIP with colourings $a_1,a_2,...,a_{s-1}$ and $a_2,a_3,...,a_s$ respectively. Suppose that for some $t$, $2\le t \le s-1$, $a_{t-1} \ne a_{t+1}$.

\noindent Let $P_1:p_1p_2,p_2p_3,...,p_{t-1}p_t$ and $P_2: p_tp_{t+1},p_{t+1}p_{t+2},...,p_{s-1}p_s$ as shown in Figure 3.

\begin{figure}[H] 
    \centering
    \includegraphics[width=15cm]{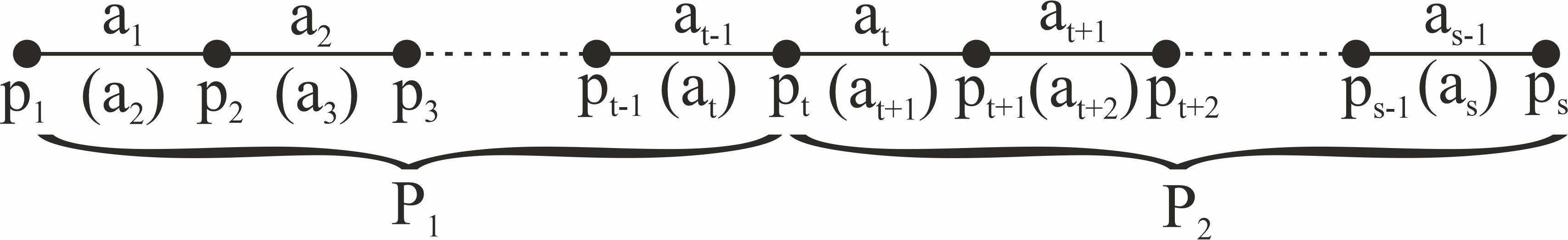}
    \caption{{\it $P$ cut into $P_1$ and $P_2$.}}
\end{figure}

\noindent Let $P_1$ have a first colouring $a_1,a_2,...,a_{t-1}$ and a second colouring $a_2,a_3,...,a_t$ respectively, and $P_2$ have the one colouring $a_{t+1},a_{t+2},...,a_s$ respectively. Since  $a_{t-1} \ne a_{t+1}$, $P_1$ is a CIP since at the vertex $P_t$ neither of the possible colourings of $p_{t-1}p_t$ clash with the colour $a_{t+1}$ assigned to $P_2$. See Figure 4.

\begin{figure}[H] 
    \centering
    \includegraphics[width=13.5cm]{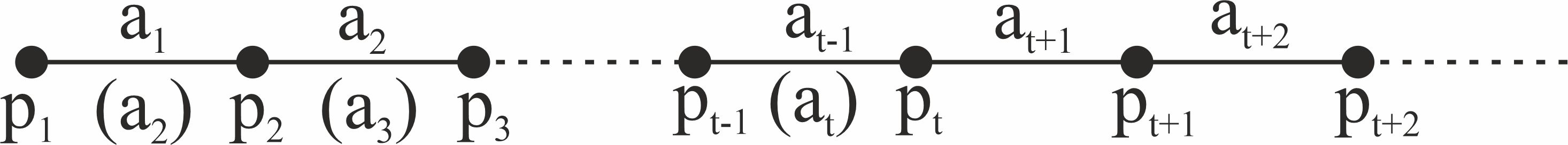}
    \caption{{\it The CIP $P_1$.}}
\end{figure}

\noindent Similarly suppose that $P_2$ has a first colouring $a_t,a_{t+1},...,a_{s-1}$ with a second colouring $a_{t+1},a_{t+2},...,a_{s}$, and that $P_1$ is coloured $a_1,a_2,...,a_{t-1}$ (with no second colouring). Since  $a_{t-1} \ne a_{t+1}$, $P_2$ is a CIP since, at the vertex $p_t$, neither of the two possible colourings of $p_tp_{t+1}$ clash with the colour  $a_{t-1}$ on the edge $p_{t-1}p_t$. (See Figure 5).
\begin{figure}[H] 
    \centering
    \includegraphics[width=14cm]{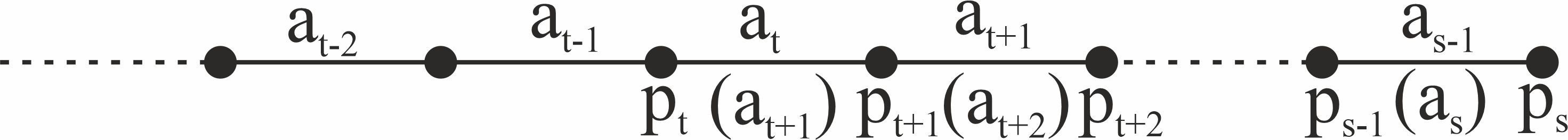}
    \caption{{\it The CIP $P_2$.}}
\end{figure}
To sum up:

\noindent {\it Lemma 9.}
\begin{enumerate}[label=(\roman*)]
\item If $a_{t-1} \ne a_{t+1}$ and if $P_2$ is coloured $a_{t+1},a_{t+2},...,a_s$ respectively, then $P_1$ is a $w$-restricted CIP with colourings $a_1,a_2,...,a_{t-1}$ and $a_2,a_3,...,a_t$ respectively.

\item If $a_{t-1} \ne a_{t+1}$ and if $P_1$ is coloured $a_1,a_2,...,a_{t-1}$ respectively, then $P_2$ is a $w$-restricted CIP with colourings  $a_t,a_{t+1},...,a_{s-1}$ and $a_{t+1},a_{t+2},...,a_{s}$ respectively.
\end{enumerate}

In Part 2  of the proof of Theorem 10 we shall extend this idea slightly by allowing $p_1=p_s$. This does not affect Lemma 9.

\section{The choice index of a simple graph}

We turn now to the proof of our main result, Theorem 5. In fact we shall prove a slightly extended version of Theorem 5.

\noindent {\it Theorem 10. Let $G$ be a finite simple graph. Let $vp_1$ be an edge of $G$. Then,}
\begin{equation*}
c'(G)\le \Delta(G)+2
\end{equation*}

\noindent {\it Moreover, given any list- assignment L to the edges of  $G$, such that  $|L(e)|\geq \Delta (G)+2$ $(\forall e \in E(G))$, then, for any vertex $w\notin \{v,p_1\}$ in $G$ and any proper L-edge-colouring of $G$, if $vp_1$ is given the colour $c$ in the proper L-edge-colouring of $G$, there is a $w$-restricted CIP starting with $vp_1$ coloured $c$}.

\noindent {\it Proof}.  First let us observe that Theorem 10 is true if $\Delta(G)=0$ or $\Delta(G)=1$. If $\Delta(G)=2$ then $G$ consists of disjoint paths and cycles. For a path $P$ we may colour the edges one by one starting at one end and it follows that $c'(p)=2$. For a cycle $C$, we may start at one edge $p_1p_2$ colour the edges one by one, with $|L(e)|=2$, on each edge, until we reach $p_1$ again, and for this find edge $e$ we may need $|L(e)|=3$. Thus $c'(G)=3$, so $c'(G)=\Delta(G)+1$. The extra requirement about the existence of a $w$-restricted CIP starting at any edge is readily seen below to be true in these cases.  Consider the case when $\Delta(G)=2$. 

If $p_1p_2$is an edge coloured $a$ and $|L(p_1p_2)|\ge4=\Delta(G)+2$, there are three colours in $L(p_1p_2)$ but not used on $p_1p_2$, and so, if $p_3$  is adjacent to $p_2$, there are two colours in $L(p_1p_2)$ not used on $p_2p_3$, so there are two colours which could be placed on $p_1p_2$ instead of $a$. So $p_1p_2$ is a CIP in its own right coloured $a$, but with a possible alternate colour $b$. The path $p_1p_2,...,p_{s-1}p_s$ reduces in this case to just the edge $p_1p_2$. The $w$-restricted requirement in this case is vacuous whichever vertex $w \notin \{p_1,p_2\}$ in $G$ is, since the CIP only has one edge.

We may note at this point that in the case when $\Delta(G)=3$, Ellingham and Goddyn \cite{MNEllingham1996} showed that $c' (G)=\chi'(G)$. We shall not use this fact, but start our general argument with the case $\Delta(G)=3$.

From now on, suppose that $\Delta(G)\ge3$. The proof is divided into two parts. Roughly, in the first part we assume that each edge can be chosen to be the start of a $w$-restricted finite colour interchange path. In the second part, we justify this assumption. We adapt the proof of  Vizing \cite{VGVizing1964} in 1964, often called the Vizing fan argument, that $\chi'(G)\le \Delta(G)+1$. The process we describe has to terminate, and so we shall describe what we do at each step, and if this does not produce the colouring we are looking for, then we go on to the next step. We suppose we have a particular simple graph G with a given maximum degree $\Delta(G)$. We shall suppose that $G$ is provided with an edge-list assignment $L$ and that $|L(e)|\ge\Delta(G)+2$ for each edge $e \in E(G)$. We shall suppose that, for some edge $f=vy$, in $E(G)$, $f$ is not coloured, but that for each $e \in E(G-f)$, $e$ is coloured with a colour in $L(e)$. We shall suppose that for each edge $e$ of $G-f$, there is a finite $w$-restricted CIP starting on $e$. We shall show that the edge colouring can be modified slightly so that $f$ can be restored to $G$ coloured with a colour in $L(f)$. Then $G$ will have an L-edge-colouring. Then we show in Part 2 that $G$ has a $w$-restricted finite CIP starting on any edge.

\noindent {\it \textbf {PART 1}}

The proof is by induction on the value of $\Delta(G)$, and for fixed value of $\Delta(G)$, this proof is by induction on the number of edges in $G$.

We showed above that the theorem is true if $\Delta(G)\leq 2$. Now suppose that $\Delta(G)\geq 3$ and that the theorem is true for all simple graphs with maximum degree at most $\Delta(G)-1$. Now let $G$ be any simple graph with $m\geq 1$ edges. We shall suppose that if H is any subgraph of $G$ with $m-1$ edges and with $\Delta(H)\leq \Delta(G)$, and L is any edge-list assignment to $H$ with $|L(e)|\geq \Delta(G)+2$ for all $ e \in E(H)$, then, for any vertex $w$ of $H$, $H$ has the $w$-restricted colour interchange property, i.e if $v$ is any vertex of $G$, $v \neq w$, and $p_2$ is any vertex adjacent to $v,p_2\neq w$, then there is a $w$-restricted colour interchange path $P,p_1p_2,p_2p_3,...,p_{s-1}p_s$  (with $v=p_1$) such that there are colours $a_1,a_2,...,a_s$ with an L-edge-colouring of $P$ with $p_ip_{i+1}$ having colour $a_i, (1\leq i\leq s-1)$ and another with $p_ip_{i+1}$ having colour $a_{i+1},(1\leq i\leq s-1)$.

At $v$ at most $\Delta(G)-1$ colours have been used, and so there are at least three colours, say $a,b,c$ in $L(vy_1)$ which have not been used. If any of those, say $a$, has not been used on $y_1$, then we may colour $f$ with $a$, and then we have an L-edge-colouring of $G$. So suppose that $a,b,c$ are all present at $y_1$. At least three colours in $L(vy_1)$ are missing at $y_1$. If  $t_1 \in L(vy_1)$ is missing at $y_1,t_1 \neq a$, and is also missing at $v$, then we may colour $f=vy$ with $t_1$, and then we have an L-edge colouring of $G$. So we may suppose that $t_1$  is missing at $y_1$  but is present at $v$, so suppose the edge $vy_2$ is coloured $t_1$. 

There are at least two colours in $L(vy_2)$ which are not used at $y_2$. Let $t_2$ be such a colour. If $t_2$ is missing at $v$ as well, then we may recolour $vy_2$  with colour $t_2$, and colour $f=vy$ with colour $t_1$. Then we have an L-edge-colouring of $G$. So suppose there is an edge $vy_3$ incident with $v$ coloured $t_3$. 

We now start on a process of construction a Vizing-type fan on $v$. We find a sequence $y_1,y_2,...,y_i$, of distinct vertices with $y_2,...,y_i$ adjacent to $v$ in $G-f$, with $vy_2,vy_3,...,vy_i$ coloured with distinct colours $t_1,...,t_{i-1}$ respectively, with $t_2,...t_i$ missing from $y_2,...,y_{i-1}$ respectively and with $t_2 \in L(v,y_1)\cap L(v,y_2),t_3 \in L(v,y_2)\cap L(v,y_3 ),...,t_{i-1} \in L(v,y_{i-1})\cap L(v,y_i)$  respectively. This is illustrated in Figure 6.

\begin{figure}[H] 
    \centering
    \includegraphics[width=9cm]{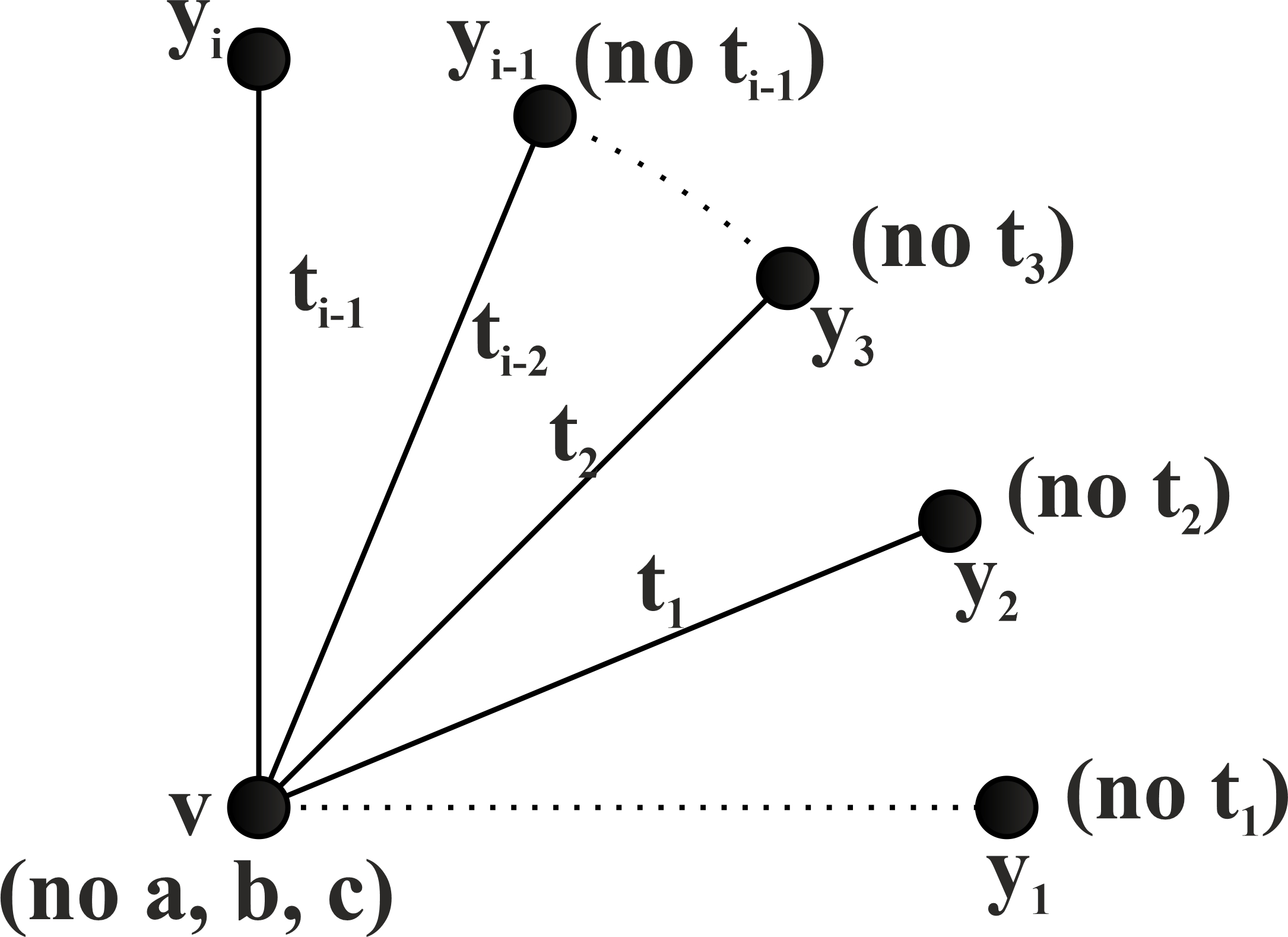}
    \caption{A Vizing fan.}
\end{figure}
 
If there is a colour $t_i \in L(v,y_i)$  which is missing at $y_i$ then either $t_i$  is missing at $v$, or $t_i$ is used on some edge incident with $v$. If $t_i$  is missing at $v$ then we recolour $vy_i,vy_{i-1},vy_{i-2},...,vy_2$ with $t_i,t_{i-1},t_{i-2},...,t_2$ respectively and colour $f=vy_1$ with $t_1$. Then $G$ has an L-edge-colouring. If $t_i$ is used on an edge $vy_{i+1}\ne vy_j$ for any $j$, $2\le j<i$, then we add the edge $vy_{i+1}$ to the fan, and continue to construct the fan.

It remains to consider the possibility that $t_i$ is used on some earlier edge of the fan (so $t_i=t_j$ for some $j$, $1\le j<i$). Since $|L(vy_i)| \ge \Delta(G) + 2$, there is some colour $a_1$ lying in $L(vy_i)$ which is absent at $v$.  If $a_1$ is absent at $y_i$ then we may obtain an L-edge-colouring of $G$ as described above. If $a_1$ is present at $y_i$, then let $p_1=y_i$ and let $p_1p_2,p_2p_3,...,p_{s-1}p_s$  be a $v$-restricted finite colour interchange path $P$ coloured  $a_1,a_2,...,a_{s-1}$ respectively for which $a_2,a_3,...,a_s$ also corresponds to an L-edge-colouring of $G-f$. Provided the following remains true {\it after} the interchange of colours on the path $P$, then we can recolour and then finish the L-edge-colouring of $G-f$. We need the following to be true:

\noindent (*) $vy_h$    $(2\le h \le i)$ has colour $t_{h-1}$ and $t_h \in L(vy_h)$ is missing at $y_h$. Also $t_1 \in L(vy_1)$ and $t_1$ is missing at $y_1$.

We can ensure that (*) is true after the interchange by taking a little care in selecting the $v$-restricted colour interchange path. The care that we take involves making use of the assumption $|L(e)|\ge \Delta(G)+2$  $(\forall e \in E(G))$. (If we had just assumed that $|L(e)|\ge \Delta (G)+1$ then we could not take take the extra care we now describe). The extra care is to ensure the following condition (X) is satisfied (here we use $v$-restriction). Here we suppose that $y_h=p_k$ (so the edge $y_hp_{k-1}$ is the same as the edge $p_kp_{k-1}$ in Figure 7 below).

\noindent (X) If $y_hp_{k-1}$ is an edge in the edge colour interchange path with colour $a_{k-1}$ prior to the interchange and $a_k \in L(y_hp_{k-1})$ is a colour missing at $y_h$, then we may only use $a_k$ on $y_hp_{k-1}$  such that $a_k\ne t_h$.

\newpage
\begin{figure}[H] 
    \centering
    \includegraphics[width=13cm]{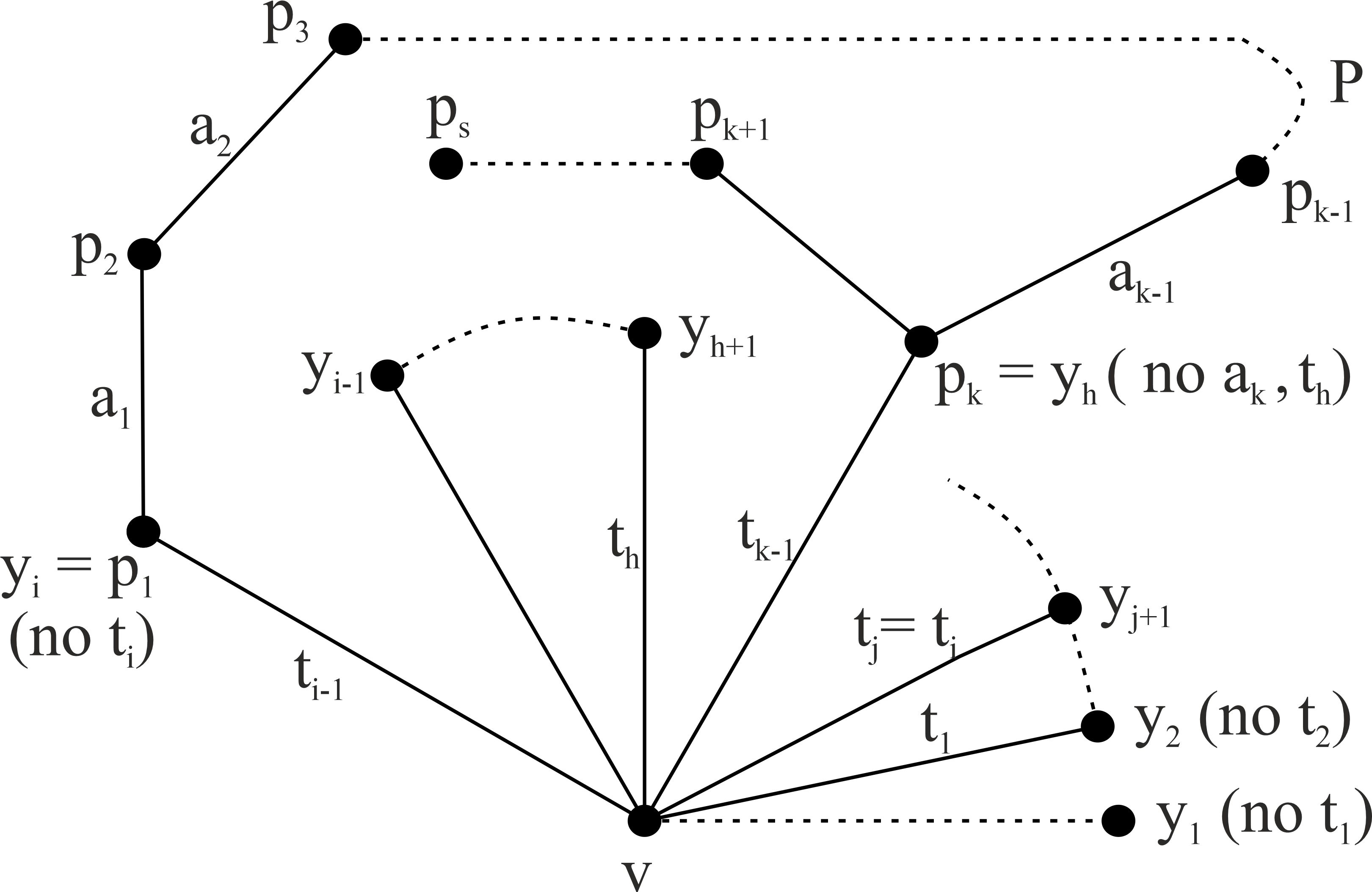}
    \caption{Condition X. Before interchanging colours on $P$.}
\end{figure}

\begin{figure}[H] 
    \centering
    \includegraphics[width=13cm]{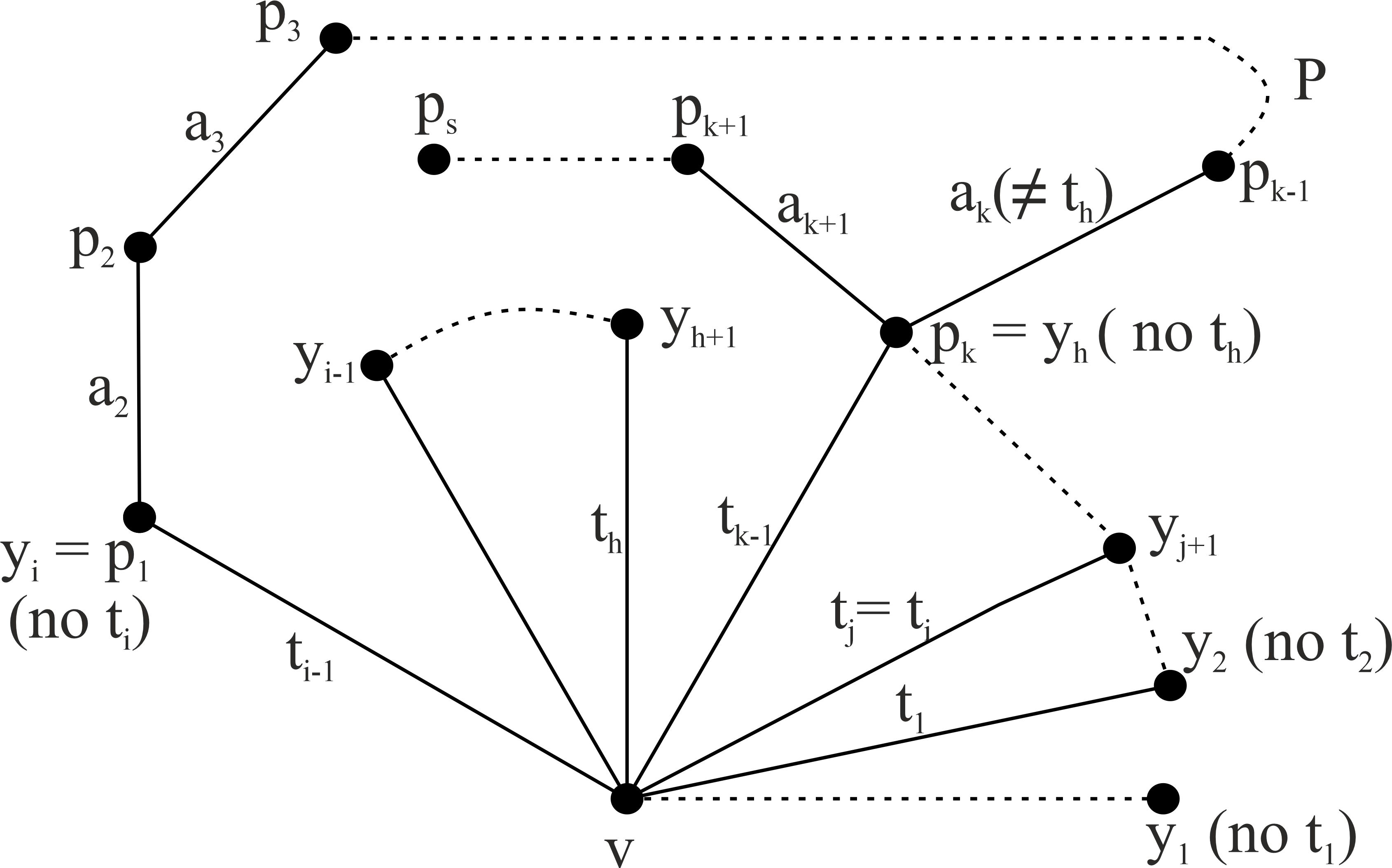}
    \caption{Condition X. After interchanging colours on $P$.}
\end{figure}

Since the degree of $y_h=p_k$ is at most $\Delta(G)$ and we are assuming that we have at least $\Delta(G)+2$  colours in each list, there is a suitable colour in $L(y_hp_{k-1})=L(p_kp_{k-1})$, as explained in Lemma 8.  Thus we can ensure that (*) is true.

So, assuming (*) is true, we may interchange the colours on the path $p_1p_2,p_2p_3,...,p_{s-1}p_s$ where $y_i=p_1$ (See Figure 8). Then each edge $p_{k-1}p_k$ $(2\le k\le s)$  receives the colour $a_k$, so $p_1p_2$ receives the colour $a_2$, and $a_1$ is no longer used on $y_i$. If $a_1$ remains missing at $v$ after the interchange then we can recolour $vy_k$ with colour $t_k  (2\le k \le i)$ and colour $f=vp_1$ with colour $t_1$. Then we have an L-edge-colouring of $G$.

We need to consider the possibility that $a_1$ is no longer missing at $v$. This might be because the edge-colour interchange path stopped at $v$ (so $p_s=v$ and $a_s=a_1$) or it might be that some internal edge of the path was incident with v and took colour $a_1$ after the interchange. See Figures 9,10 and 11 which illustrate these two possibilities.

Let us consider the first possibility illustrated in Figure 9. In this case the situation we are faced with is depicted in the “before” part of Figure 9. We interchange on the path, after which $a_1$ is missing at $y_i$, $a_{s-1}$  is missing at $v$, but $a_s=a_1$ is present at $v$. Then we rename $p_{s-1}$ as $y_{i+1}$. Now we rename $a_1$ as $t_i$, so that $t_i$ is present at $v$ and we use it to colour the edge $vy_{i+1}$. Note that, by our original argument, there are still at least three colours in $L(vy_1)$ which are missing at $v$, even if $a_1$ is no longer one of them.

\begin{figure}[H]
  \centering
  \begin{tabular}{@{}c@{}}
    \includegraphics[width=11cm]{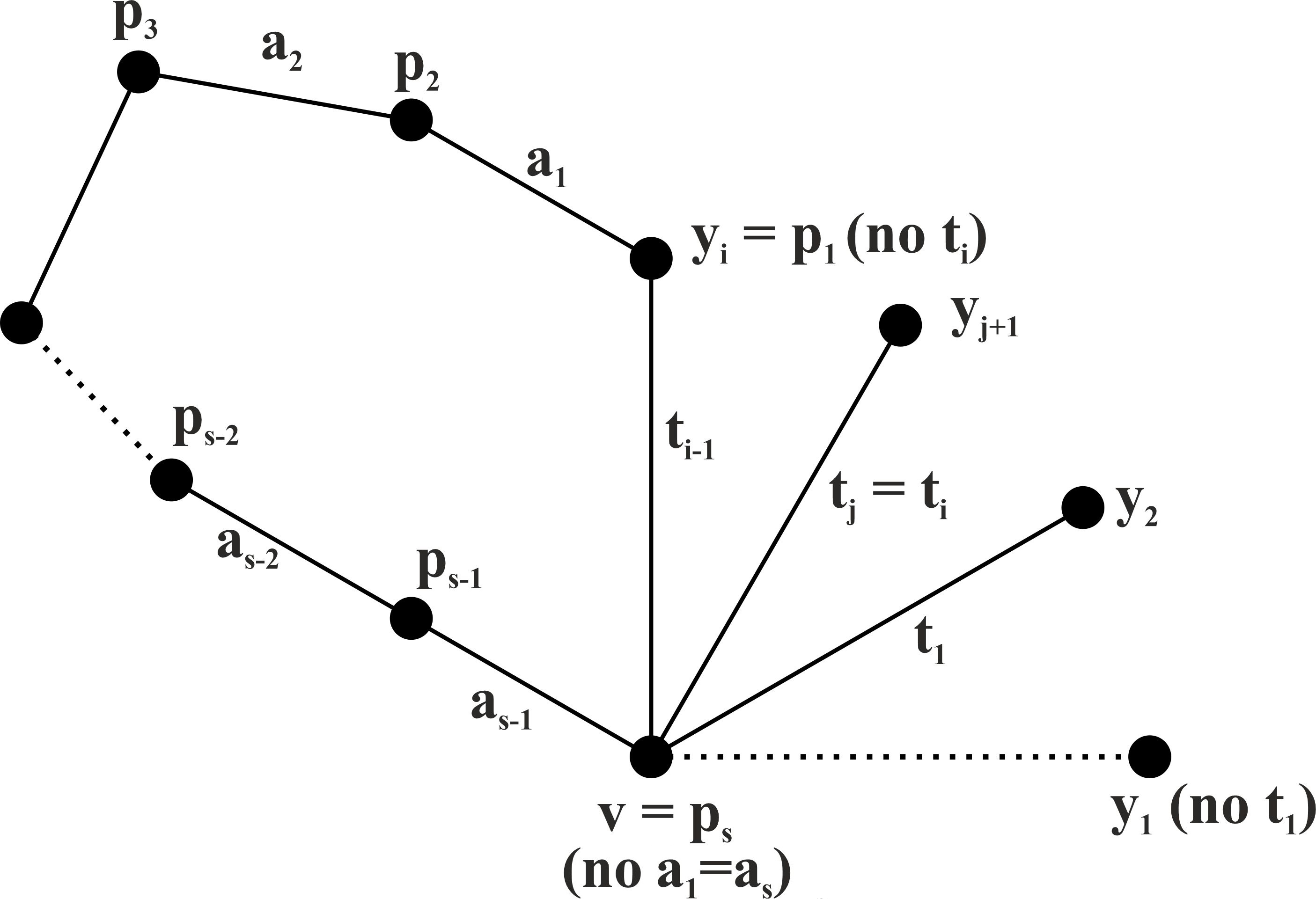} \\[\abovecaptionskip]
    Before
  \end{tabular}
  \begin{tabular}{@{}c@{}}
    \includegraphics[width=11cm]{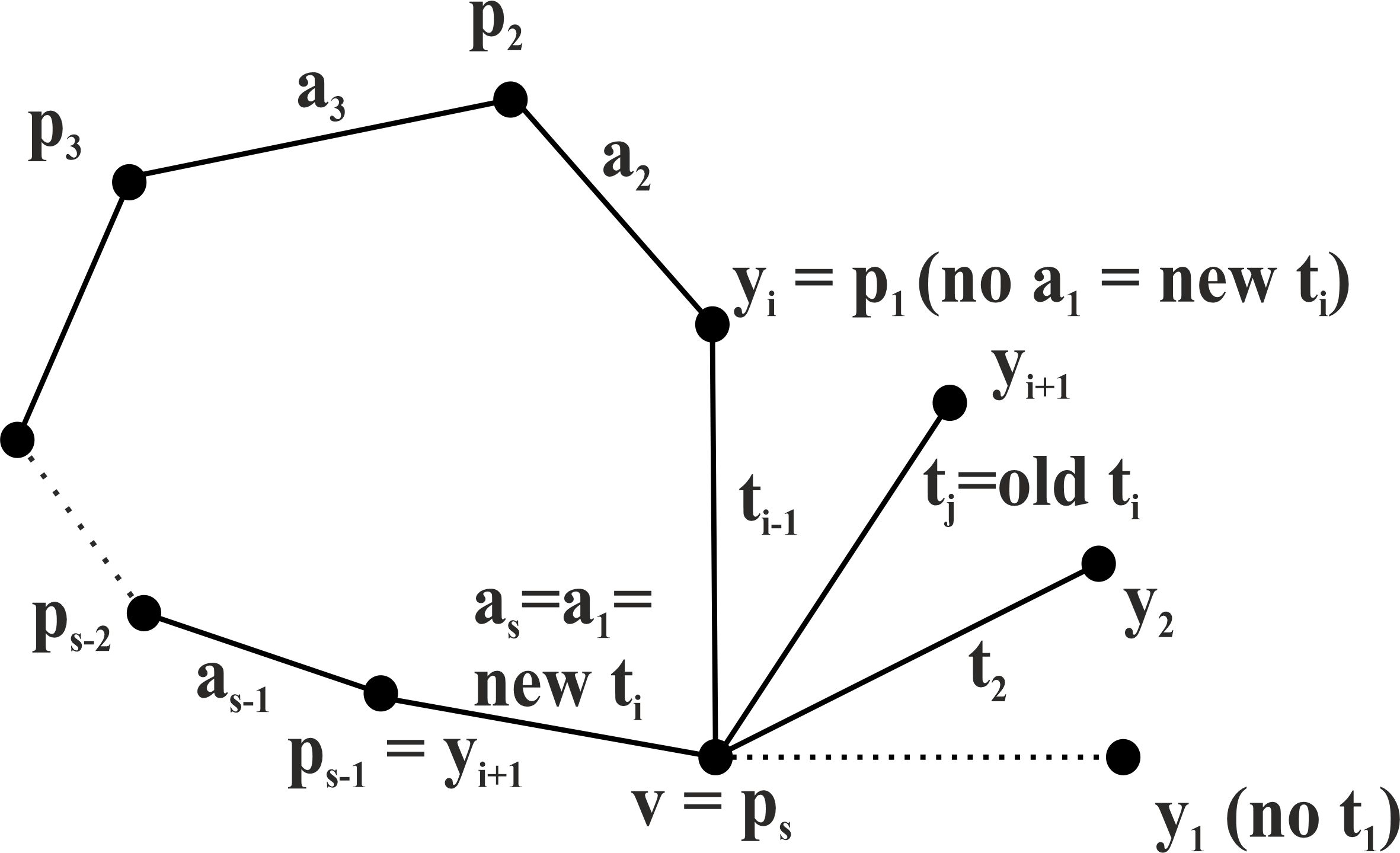} \\[\abovecaptionskip]
   After
  \end{tabular}

  \caption{Before and After when $p_s=v.$}\label{fig:myfig}
\end{figure}

Now consider the possibility that some internal edge of the path was incident with $v$ and takes colour $a_1$ after the interchange. This can happen in just one way (note that we have excluded the possibility that $a_1=a_r$). The only possibility is that $a_1= a_{r+1}$.

\begin{figure}[H] 
    \centering
    \includegraphics[width=11.5cm]{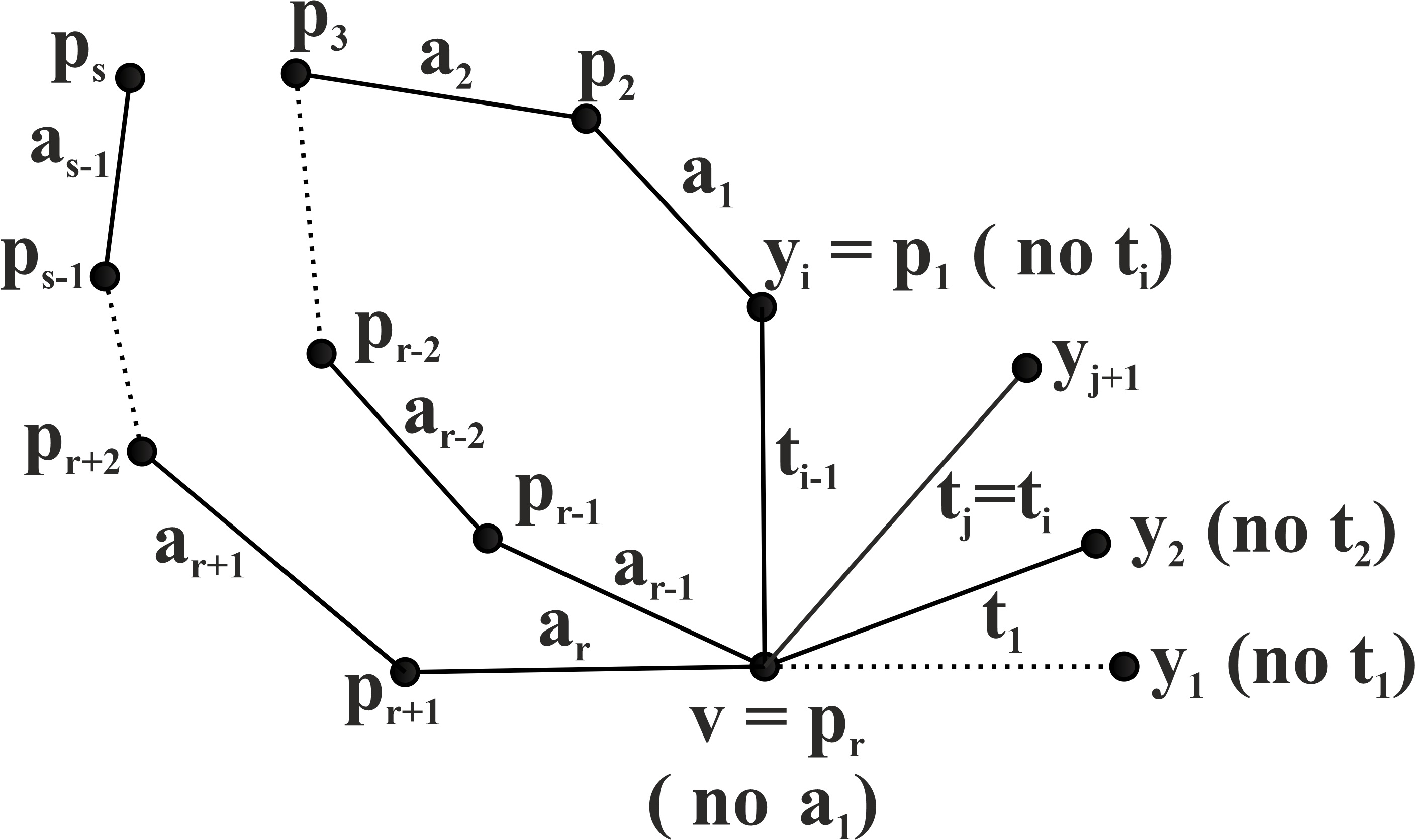}
    \caption{$v$ is an internal vertex of the interchange path.}
\end{figure}

\noindent Interchange colours on the whole path starting on $p_1p_2$. This is illustrated in Figure 11.

\begin{figure}[H] 
\centering
\includegraphics[width=11cm]{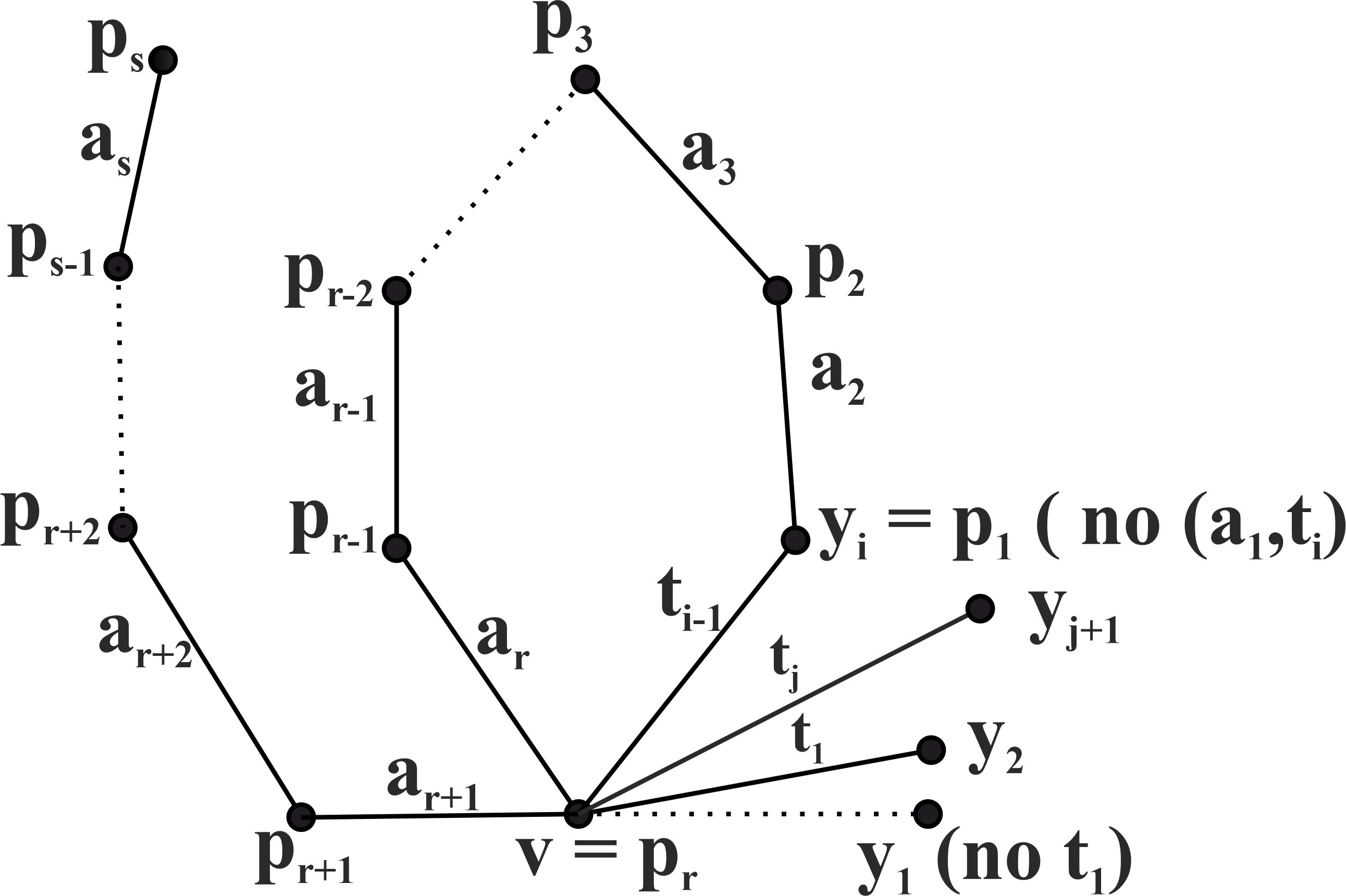}
\caption{The colour interchange.}
\end{figure}

In this case, rename  $p_{r+1}$ as $y_{i+1}$ and rename $a_{r+1}=a_1$ as $t_i$. This new $t_i$ is missing at $y_i$ and colours $vy_{i+1}$.

It is worth pointing out at this point that $a_1$ need not be the original $a$ that we identified at the start of the recolouring process. All that is needed for the original assignment to work is that $a_1$ is present at $y_i$ and is not present at $v$. Then it will follow, after the various changes made in the proof, that $t_i \notin \{t_1,…,t_{i-1}\}$ and $t_i \ne t_j$. 

Thus at this point, provided we have adhered to (*) in any path used, either we can colour $vy_1$ or we can find an L-colouring of $G$ in which $vy_1$ is uncoloured and $vy_i$ has colour $t_{i-1}$ and $t_i \in L(vy_i)$  is missing at $y_i$ for all $i$, $1\leq i \leq d$, where $d=d_G(v)$. It remains to show that in this case also we can colour $vy_1$.

At least two colours in $L(vy_d)$ are missing at $v$. Let $b_1$ be such a colour. Let $y_d$ be labelled $q_1$ and $q_1q_2,q_2q_3,...,q_{x-1}q_x$ be a $v$-restricted path coloured $b_1,b_2,...,b_{x-1}$ respectively have the property that $b_2,b_3,...,b_x$ also corresponds to an L-edge-colouring of $G$. We may assume that this path is constructed so that (X) is satisfied. Unless this path assigns $b_1$ to an edge incident with $v$, then we can interchange colours on the path, giving $vy_d$ colour $b_1$, and we can recolour $vy_{d-1}$ with colour $t_{d-1}$, $vy_{d-2}$ colour $t_{d-2},...,vy_2$ colour $t_2$ and colour  $vy_1$ with colour $t_1$. This is then an L-edge-colouring of $G$.

Next suppose that $v=q_{x}$ and that $b_1=b_x$, which is one way in which $b_{1}$ might be assigned to an edge incident with $v$ (See Figure 12).

\begin{figure}[H] 
\centering
\includegraphics[width=11cm]{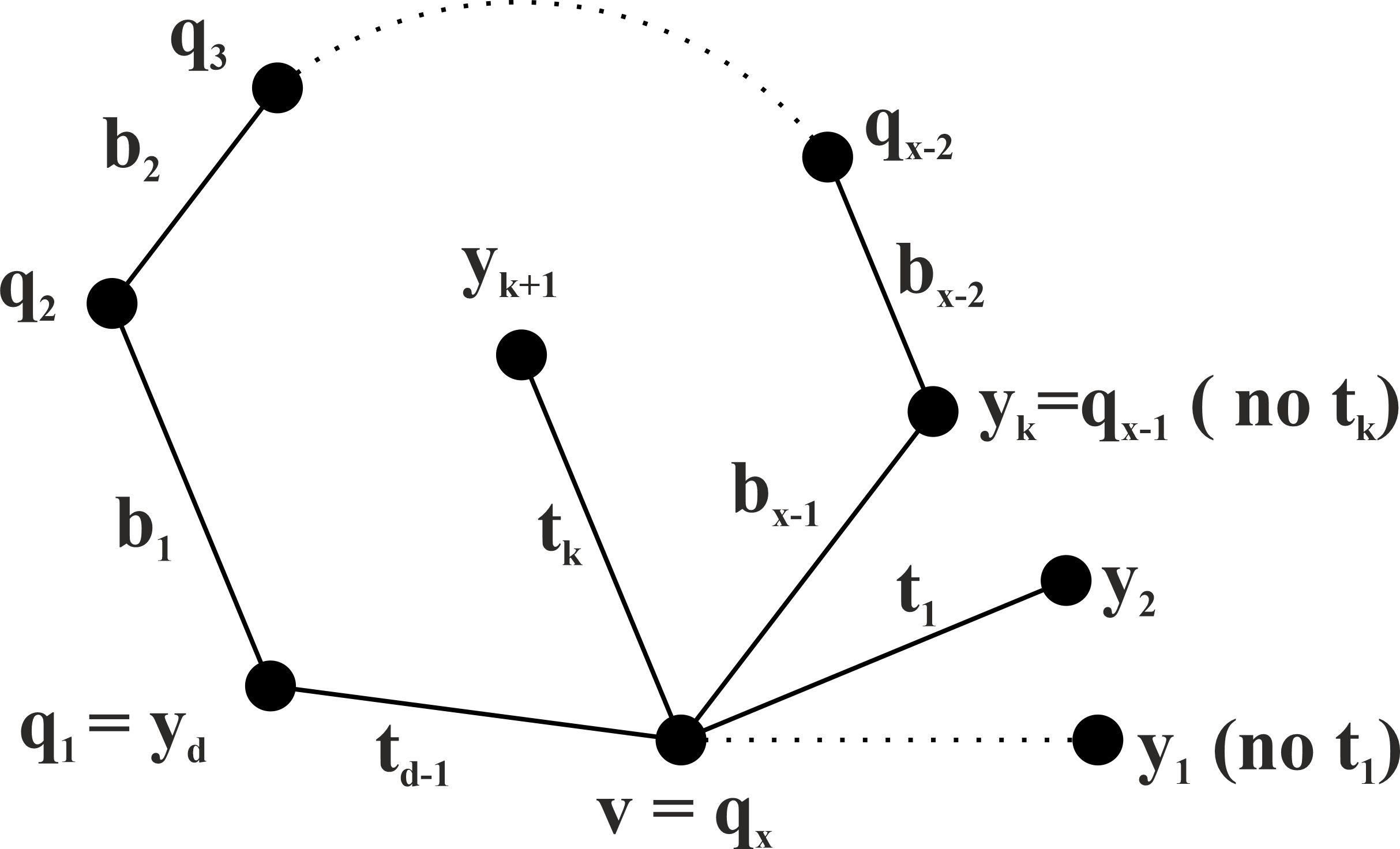}
\caption{$v= q_x$ before interchange on the path.}
\end{figure}

We interchange colours on the path $q_1q_2,q_2q_3,...,q_{x-1}q_x$. Then edge $vy_k=q_{x-1}q_x$ which was coloured $b_{x-1}=t_{k-1}$  is now coloured $b_1$ (See Figure 13). Then we recolour $vy_{k-1}$ with $t_{k-1}$, $vy_{k-2}$ with $t_{k-2},...,vy_2$ with $t_2$ and $vy_1$ with $t_1$. This gives an L-edge-colouring of $G$.

\begin{figure}[H] 
\centering
\includegraphics[width=10cm]{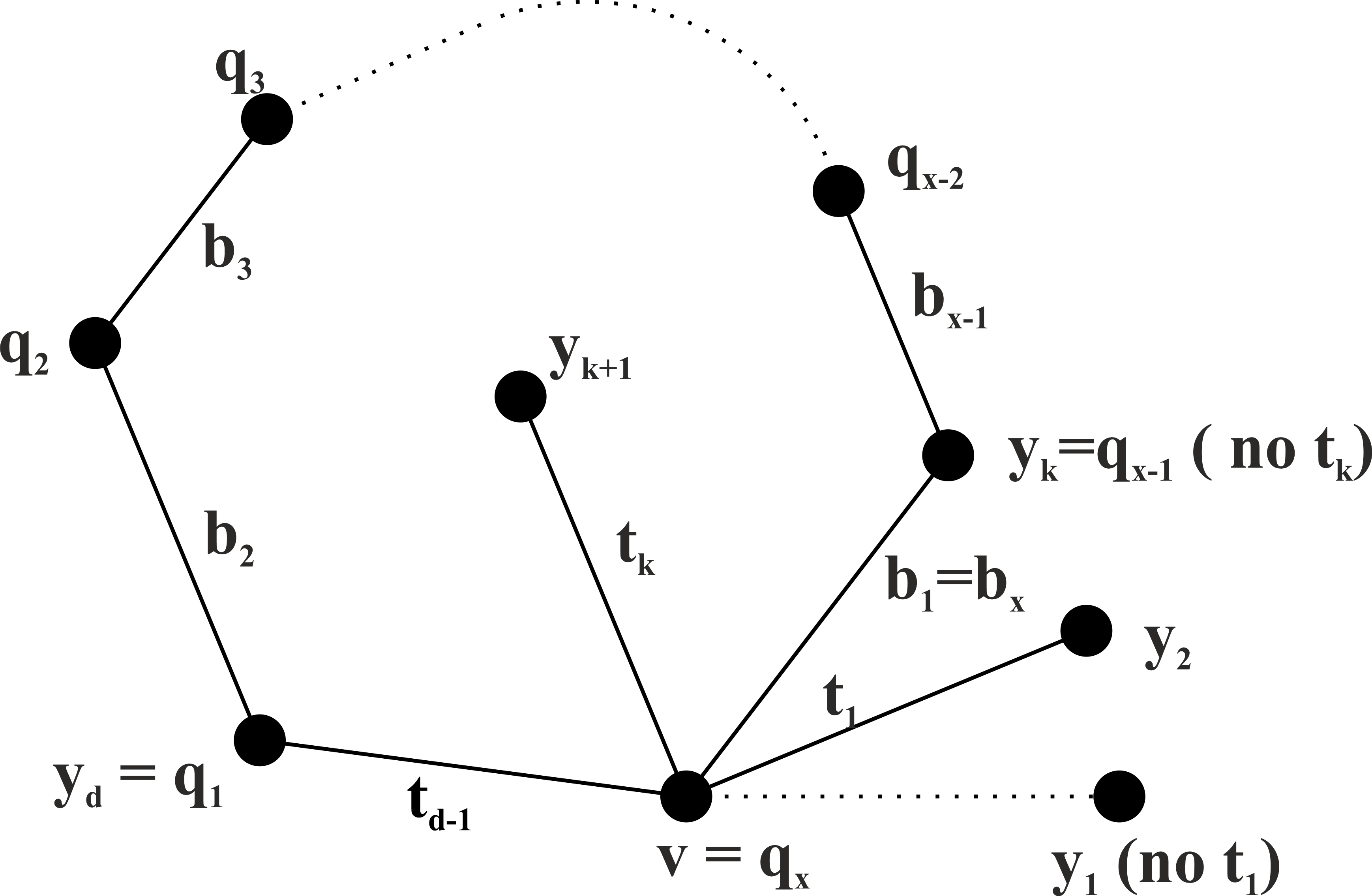}
\caption{After interchanging on the path when $v=q_x$.}
\end{figure}

The other way that $b_1$ could be assigned to an edge incident with $v$ is if $v$ is incident with an internal edge of the path. See Figure 14 to illustrate the situation.We can suppose

\begin{figure}[H] 
\centering
\includegraphics[width=13cm]{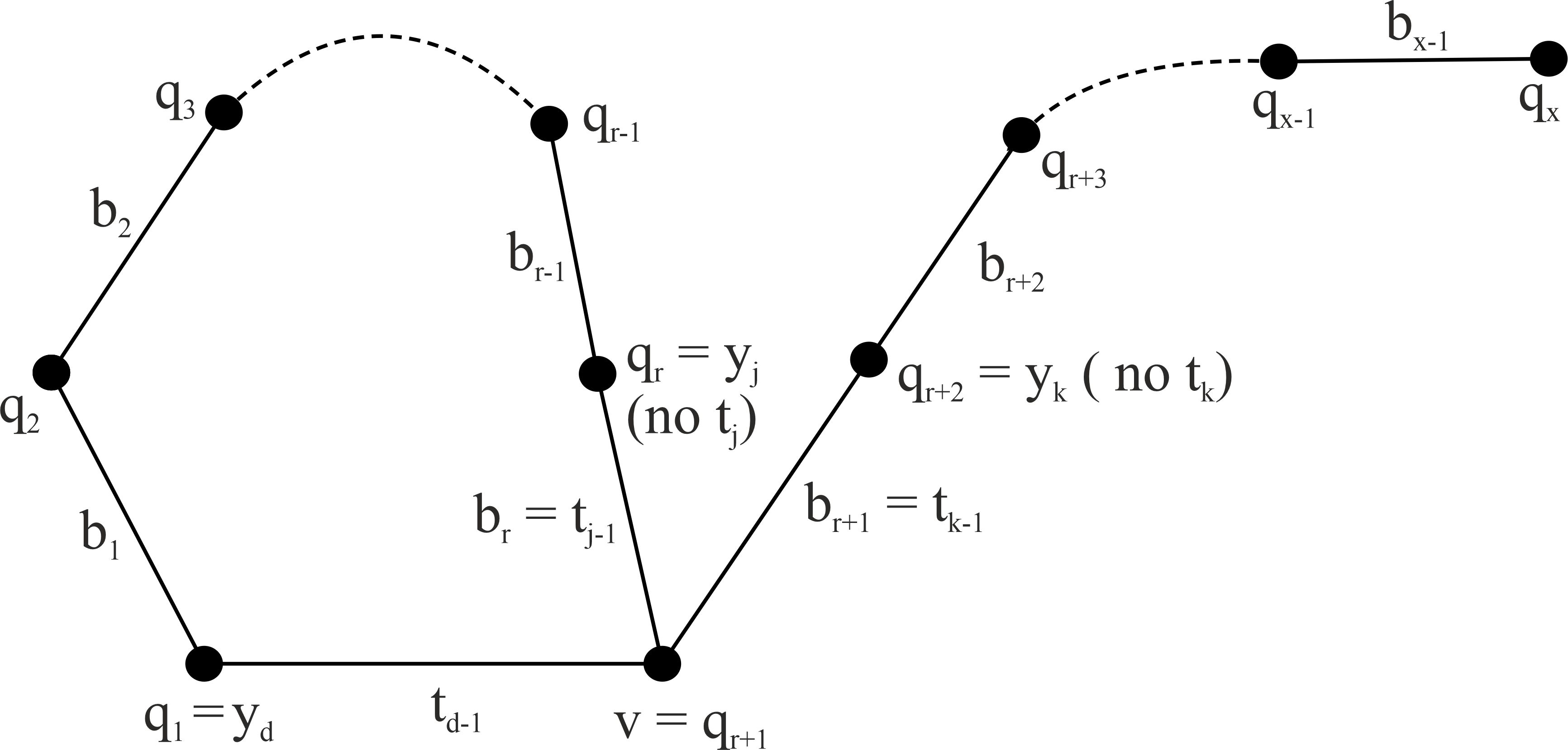}
\caption{$v=q_ {r+1}$ and $q_{r+1}=y_k$  in the interchange path.}
\end{figure}

\noindent that $v=q_{r+1}$ and $vy_k=vq_{r+2}$. We may call the first part of the path  $q_1q_2,q_2 q_3,...,q_{r-1}q_r$ the `head' of the path, and the second part $q_{r+2}q_{r+3},q_{r+3}q_{r+4},...,q_{x-1}q_{x}$ the `tail'. Since $b_1$ is missing at $v=q_{r+1}$, we know that $b_1\notin \{b_r, b_{r+1}\}$. If we were to recolour the whole path, then the head would be recoloured $b_2,b_3,...,b_{r}$ and the tail would be recoloured $b_{r+3},b_{r+4},...,b_x$. The central part, $q_rq_{r+1},q_{r+1}q_{r+2}$ would be recoloured $b_{r+1},b_{r+2}$. Therefore the only way $b_1$ could occur on an edge incident with $v$ after the recolouring is if $b_1=b_{r+2}$. The edge $q_{r+1}q_{r+2}$ is coloured $t_{k-1}$ before the recolouring, so $t_{k-1}=b_{r+1}$, and we may suppose that the edge $vq_r$ is coloured $t_{j-1}$ before the recolouring. We may suppose that $q_{r+2}=y_k$ and has no edge coloured $t_k$ incident with it, and that $q_r=y_j$ and has no edge coloured $t_j$ incident with it. How we proceed at this point depends upon whether $j<k$ or $j>k$. 

Suppose that $k<j$. We interchange the colours on the path $q_{r+1}q_{r+2},q_{r+2}q_{r+3},...,q_{x-1}q_{x}$. These edge are then recoloured $b_{r+2}=b_1, b_{r+3},...,b_{x-1}, b_x$ respectively. Then the colour $b_{r+1}=t_{k-1}$ is missing at $v=q_{r+1}$ since, in particular, $t_{k-1}\ne b_{r+2}=b_1$ and $t_{k-1}\ne b_r$. This colouring is illustrated in Figure 15.

\begin{figure}[H] 
\centering
\includegraphics[width=13.5cm]{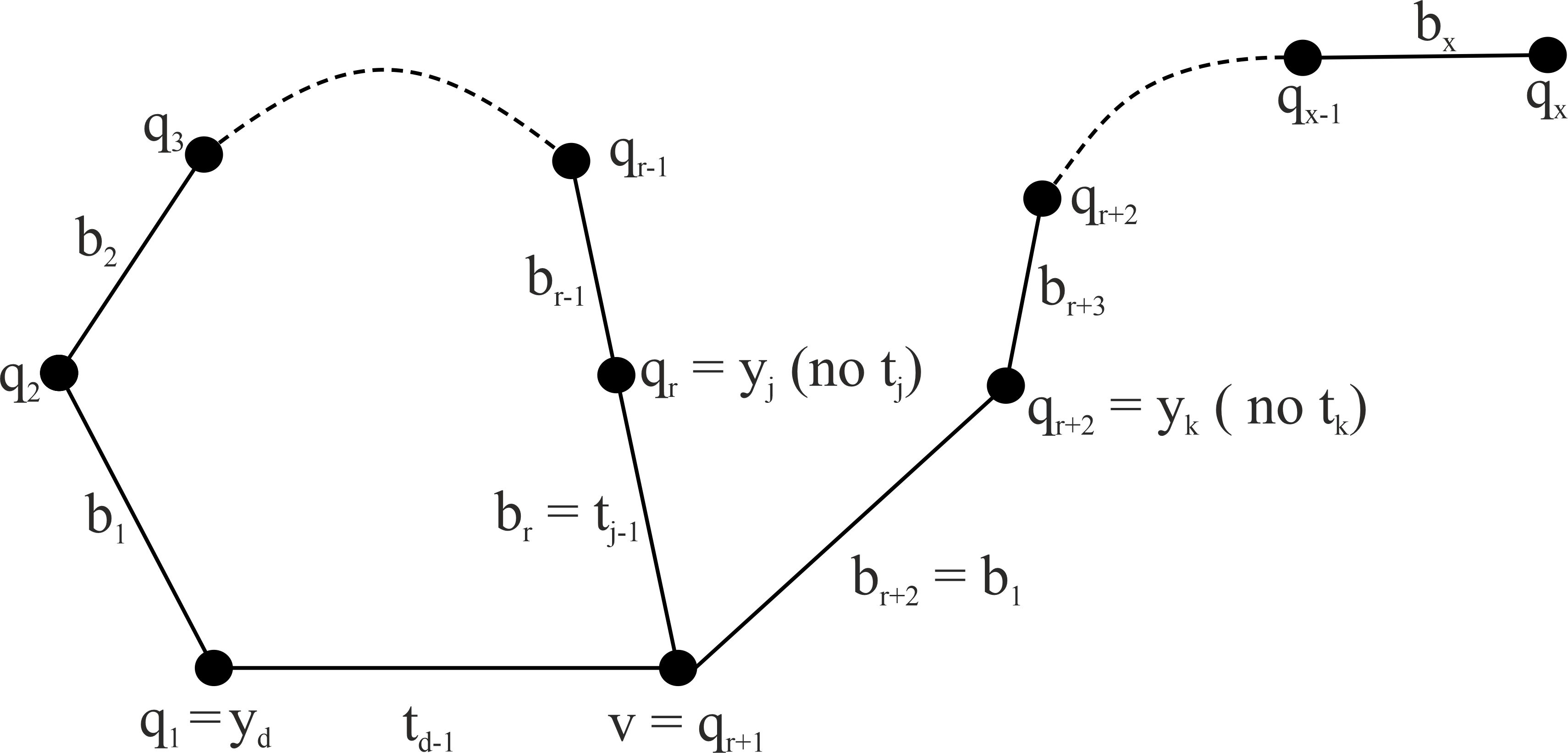}
\caption{$v=q_ {r+1}$, the colour on $q_{r+1}q_{r+2},q_{r+2}q_{r+3},...,q_{x-1}q_{x}$ are interchanged. No $t_{k-1}$ on $v=q_ {r+1}$.}
\end{figure}

We complete the L-edge-colouring of $G$ by recolouring $vy_{k-1}$ with $t_{k-1},vy_{k-2}$ with $t_{k-2},...,vy_1$ with $t_1$.

Finally suppose that $j<k$. Recall that $b_{r+2}=b_1$ and no edge coloured $b_1$ is incident with $v$. Therefore $b_{r+2} \ne t_{k-1}$. We interchange the colours on the whole path $q_1q_2,q_2 q_3,...,q_{x-1}q_x$. The situation is illustrated in Figure 16. 

\begin{figure}[H] 
\centering
\includegraphics[width=12.5cm]{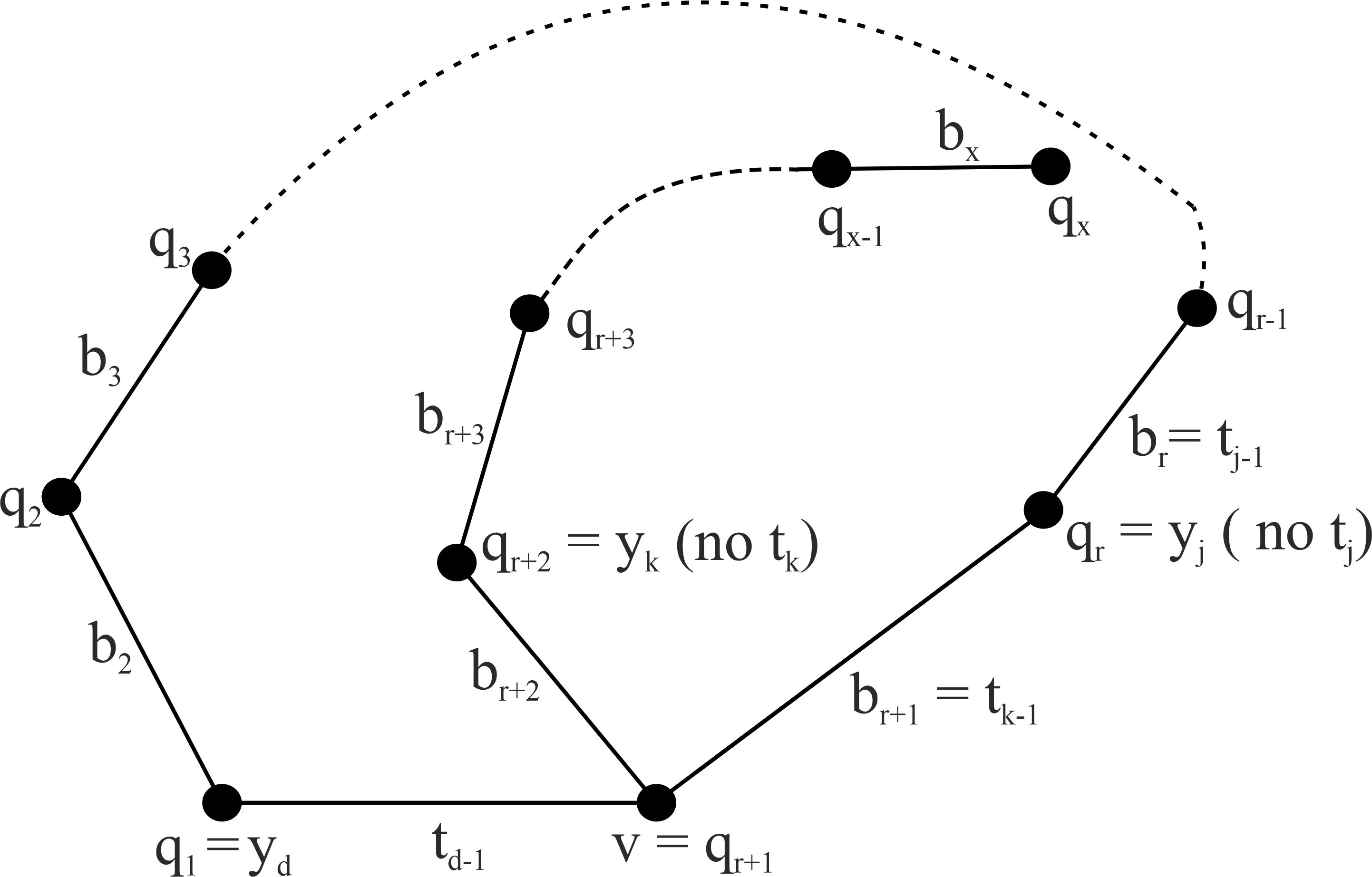}
\caption{$v=q_ {r+1}$. The colours on the whole path are interchanged. There is no $t_{j-1}$ on $v=q_ {r+1}$.}
\end{figure}
There is no edge incident with $v=q_{r+1}$ coloured $t_{j-1}$. We complete the L-edge-colouring of $G$ by recolouring $vy_{j-1}$ with $t_{j-1},vy_{j-2}$ with $t_{j-2},...,vy_1$ with $t_1$. 

At this point we have completed the induction step and shown that, assuming that $c'(G-e)\le \Delta(G)+2$ and that $L$ is an edge list assignment to $G$ and that starting with any edge $p_1p_2$, there is a $v$-restricted colour interchange path $p_1 p_2,p_2p_3,...,p_{s-1}p_s$, it follows that the L-edge colouring can be extended to include $e$ so that we have an L-edge-colouring of G. We now need to provide the second part of our proof.

\noindent {\it \textbf{PART 2}}

In this part, we need to show that, if $G$ has an L-edge-assignment such that $\Delta(G)+2 \leq |L(e)|$ $(\forall e \in E(G))$, if  $v \in V(G)$ and if $G$ has a proper L-edge-colouring, and if $xp_1$ is an arbitrary edge of $G$ with $x\neq v$, $p_1\neq v$, then $G$ has a $v$-restricted CIP starting on the edge $xp_1$ coloured $a$  if $xp_1$ receives the colour $a$ in the L-edge-colouring of $G$.

As an induction hypothesis, we can assume that $v \in V(G)$, that $G$ has a proper L-edge-colouring (by Part I), and that for any edges $e$, $f$ of $G$ with $v$ not incident with $e$ or $f$, $G-e$ has a $v$-restricted CIP starting on $f$ with $f$ receiving the colour it receives in the L-edge colouring of $G$. The induction step consists of showing that $G$ itself has a $v$-restricted CIP starting with the edge $xp_1$ coloured with the colour it received in the L-edge-colouring of $G$.

Let $y \in V(G)$ and let $xy$ be an edge of $G$. Let $a$ be a colour in $L(xp_1)$ which is absent at $x$. Let $e=xy$ and let $f=xp_1$. By induction $G-\{xy\}$ has a $v$-restricted CIP $P: xp_1,p_1p_2,p_2p_3,...,p_{s-1}p_s$ with $xp_1$ receiving the colour it received in the L-edge-colouring of $G$ (say $a$). Let $xp_1,p_1p_2,...,p_{s-1}p_s$ be coloured $a,b_1,b_2,...,b_{s-1}$, with an alternative L-edge-colouring $b_1,b_2,...,b_s$, respectively. Then $b_1$ is also absent at $x$ and $b_1 \in L(xp_1)$.

Now consider the graph $G-\{xp_1\}$. Let $e=xp_1$ and let $p_1=r_1$. By induction in $G-\{xp_1\}$ there is a $v$-restricted CIP $R: r_1r_2,r_2r_3,...,r_{u-1}r_u$ starting on the edge $r_1r_2$ coloured $b_1$. The path $R$ might or might not be equl to the path $P$ (with edge $xp_1$ removed)- it does not matter. We may suppose that $R$ is coloured $d_1,d_2,...,d_{u-1}$ or $d_2,d_3,....,d_u$ respectively. If $R$ does not pass through, or end on $x$, then we can add the edge $xr_1$ to the start of $R$ with colour $a$ and obtain a $v$-restricted CIP in $G$ itself with L-edge-colouring $a,d_1,d_2,...,d_{u-1}$ and an alternative colouring $d_1,d_2,...,d_u$.

But suppose that $R$ does pass through, or end on $x$. Then consider the path $R$ with the extra edge $xr_1$ at the start. Call this $R^{+}$. $R^{+}$ cannot be a CIP since $R^{+}$ is not a path. But we can `cut' it along the lines explained in Section 3 at the vertex  $r_1$ provided $a \ne d_2$. Thus if $a \ne d_2$ we obtain

\begin{figure}[H] 
\centering
\includegraphics[width=14cm]{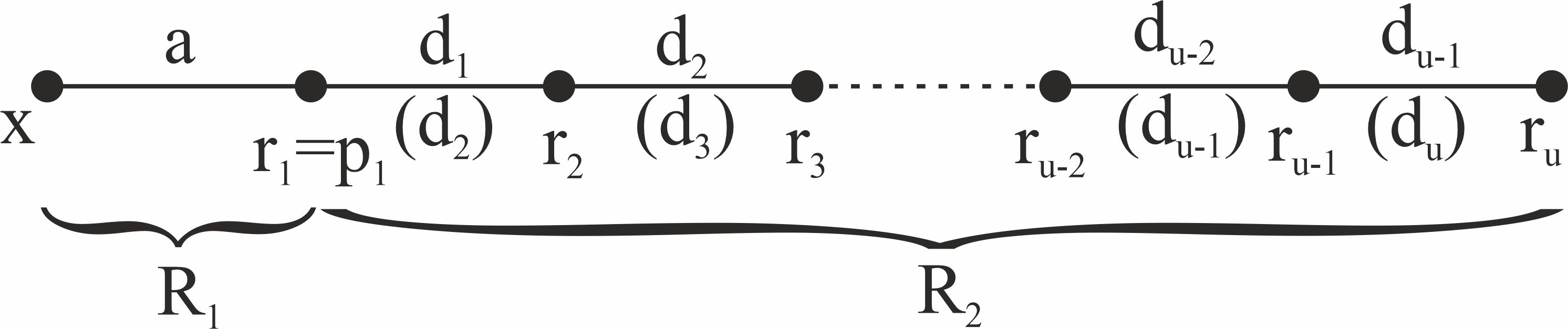}
\caption{$R^{+}$ cut at the vertex $r_1$ giving $R_1$ and $R_2$.}
\end{figure}

\noindent a $v$-restricted CIP $R_1:xp_1$ in $G$ coloured with a first colouring $a$ and a second colouring $d_1$. We also have $R_2$ coloured $d_2,d_3,...,d_u$ (See Figure 17).

We cannot do this if $d_2=a$. Let $x=r_w$. It is not the case that the colours used on $r_{w-1}r_w$ and $r_wr_{w+1}$ include $a$ and $b_1=d_1$ since $a$ and $b_1$ are absent at $x$. If there were no index $t$ such that $d_{t-1} \ne d_{t+1}$, then $R:  r_1r_2, r_2r_3, ... ,r_{u-1}r_u$ would be coloured $b_1, a, b_1, a .b_1$, ... so both $a$ and $b_1$ would be used on $x$, which is impossible. So there is an index $t$ such that $d_{t-1} \ne d_{t+1}$. In that case we have a $v$-restricted CIP $R_{1}^{`}: xr_1,r_1r_2,...,r_{t-1}r_t$ with a first colouring $a,d_1,d_2,...,d_{t-1}$ and a second colouring $d_1,d_2,...,d_t$, if we colour $R_{2}^{`}: r_tr_{t+1},r_{t+1}r_{t+2},...,r_{u-1}r_u$ with $d_{t+1},d_{t+2},...,d_u$ respectively. Thus again we have a $v$-restricted CIP in $G$ with first edge $xp_1=xr_1$ coloured $a$. 

This finishes Part II, and Theorem 10 now follows by induction.

\section{Concluding Remarks about the choice index of a simple graph.}

The obvious question raised by our theorem and its proof is whether Theorem 10 can be improved to “$c'(G)\le \Delta(G)+1$ whenever $G$ is simple graph”. We use the extra freedom, which arises when we assume $c'(G)\le \Delta(G)+2$ in several places, but it is not obvious that the same kind of proof could not be found if we assume $c'(G)\le \Delta(G)+1$. But the difficulties which seem to arise in attempting this, also raise the question of whether there are simple graphs for which $c' (G)=\Delta(G)+2$.

It might be that our proof could be shortened, but one reason for sticking to the present proof is that it shows up quite well where the difficulties lie in trying to improve it so as to show that $c'(G)\le \Delta(G)+1$, and we hope that the present proof will be helpful in obtaining this improved bound.

To prove $c'(G)=\chi'(G)$, a different method of proof would have to be found.

\section{The total chromatic number of a simple graph.}

The total chromatic number $\chi_T(G)$ of a simple graph $G$ is the least number $j$ of colours needed to colour the edges and vertices of a graph $G$ so that 

\begin{enumerate}[label=(\roman*)]
\item no colour is used on two edges which have a common vertex;
\item no colour is used on two adjacent vertices;
\item no colour is used on a vertex $v$ and an edge incident with $v$.
\end{enumerate}

\noindent Clearly $\Delta(G)+1\le \chi_T(G)$. In 1965 a conjecture, due independently to Behzad \cite{MBehzad1965Graphs} and to Vizing \cite{VGVizing1965} was made that, if G is a simple graph, then $\chi_T(G)\le \Delta(G)+2$. For a full discussion of the attribution of this conjecture, see [\cite{ASofier}, section 16.2]. In 1993, Hilton and Hind \cite{AJWHilton1993} showed that $\chi_T(G) \leq \Delta(G)+2$, if $\Delta(G) \geq \frac{3}{4}|V(G)|$. Some related results were given in 1991 by Chetwynd, Hilton and Zhao Cheng in \cite{Chetwynd1991The}. It has been known since 1998 that $\chi_T(G)\le \Delta(G)+c$, where $c$ is a constant: Molloy and Reed \cite{MMolloy1998} showed that for large enough $\Delta(G),\chi_T(G)\le \Delta(G)+10^{26}$. Here we prove the following similar result which is slightly weaker than the conjectured result.

\noindent {\it Theorem 11. Let $G$  be a simple graph, then $\chi_T (G)\le \Delta(G)+4$}

It has been known ever since the choice index was thought of, that a bound for the choice index yields a bound for the total chromatic number. We use this short argument to prove Theorem 11.

\noindent {\it Proof of Theorem 11}. Let $\zeta$  be set of $\Delta(G)+4$ colours. Colour the vertices of $G$ with colours from $\zeta$. [It is clear by a greedy algorithm that $\Delta(G)+1$ colours suffice to colour $V(G)$]. For each edge $e$ of $G$, let $L(e)$ be a list assignment consisting of $\zeta$ less the two colours used to colour the end-vertices of $e$. Then $|L(e)|= \Delta(G)+2$  $(\forall e \in E(G))$. By Theorem 5, there is an edge-list colouring of $E(G)$. The edge-list colouring together with the vertex colouring together constitute a total colouring of $G$ with at most $\Delta(G)+4$ colours.

\section{The total list chromatic number of a simple graph}

We show in this section that virtually the same argument will prove the same bound for the list analogue of the total chromatic number.

A {\it total list assignment} to a simple graph $G$, or list assignment to $E(G)\cup V(G)$ is a function $E(G)\cup V(G)\to2^{\zeta}$. If $\Lambda$ is a list assignment to $E(G)\cup V(G)$, a proper $\Lambda$-colouring of $E(G)\cup V(G)$,  is a function $\Psi:E(G)\cup V(G)\to \zeta$,  satisfying

\begin{enumerate}
\item $\Psi(e)\in \Lambda(e)\hspace{5pt}(\forall e \in E(G))$
and $\Psi(v)\in \Lambda(v)\hspace{5pt}(\forall v \in V(G))$;
\item If $e$, $f \in E(G)$ and $e$ and $f$ have a vertex in common, then $\Psi(e)\ne \Psi(f)$;
\item If $v$, $w \in V(G)$ and $v$ and $w$ are adjacent, then  $\Psi(v)\ne \Psi(w)$;
\item If  $v \in V(G)$ and $e \in E(G)$ and $e$ is incident with $v$, then $\Psi(v) \ne \Psi(e)$.
\end{enumerate}

The {\it total choice number or total list chromatic number,} $c_T (G)$ is least number $n_1$, such that, wherever $\Lambda$ is a total list assignment to $E(G)\cup V(G)$ with $|\Lambda(e)|\ge n_1$    $(\forall e \in E(G))$ and $|\Lambda(e)|\ge n_1$   $(\forall v \in V(G))$, then there exists a proper $\Lambda$ - colouring of $E(G)\cup V(G)$.

In the case when $\Lambda(e)=\Lambda(f)=\Lambda(v)=\Lambda(w)$ for all $e$, $f \in E(G)$ and $v, w \in V(G)$ ( so the lists are all the same), then $n_1$  is the total chromatic number of $G$, and is denoted by $\chi_T(G)$.

A natural analogue of the edge-list colouring conjecture, Conjecture 2, is: 

\noindent {\it Conjecture 12. For a simple graph} 
\begin{equation*}
c_T(G)=\chi_T(G).
\end{equation*}

\noindent This conjecture seems to have been made in about 1999 by Borodin, Kostochka and Woodall \cite{Borodin1997List}, Juvan, Mohar and Skrekovski \cite{MJuvan1994} and Hilton and Johnson \cite{AJWHilton1999}.

\noindent Of course, there is really no known reason why this should not hold for multigraphs also.

\noindent We prove here:

\noindent {\it Theorem 13. For a simple graph $G$}
\begin{equation*}
\Delta(G)+1\le \chi_T(G)\le c_T(G) \le \Delta(G)+4.
\end{equation*}

\noindent {\it Proof}. The first two inequalities are obvious. To prove the third inequality. Let $\Lambda$ be a list assignment to $E(G)\cup V(G)$ with $|\Lambda(e)|\ge \Delta(G)+4$ 
$(\forall e\in E(G))$ and $|\Lambda(v)|\ge\Delta(G)+4$
$(\forall v\in V(G))$. We may properly colour the vertices of $G$ using colours from the lists $\Lambda(v)$. Since $|\Lambda(v)|\ge \Delta(G)+4$, this is clearly possible (by an obvious greedy algorithm, this is possible if $|\Lambda(v)|\ge \Delta(G)+1$ 
$(\forall v \in V(G))$. For $e \in E(G)$, let $\Lambda^*(e)$ be $\Lambda(e)$ less the two colours selected to properly colour the end vertices. Then $|\Lambda^*(e)|\ge \Delta(G)+2$ $(\forall e \in E(G))$. By Theorem 5, there is an edge-list colouring of $E(G)$. This edge-list colouring together with the vertex-colouring together constitute a total list colouring of $G$.

\section{The Hall index and the Hall condition index of a simple graph.}

We should remark that a large part of this section holds {\it mutatis mutandis} for graphs in general or for the edge-sets of multigraphs (see Hilton and Johnson \cite{AJWHilton1999}).The connection between Hall’s theorem and list colouring was noticed in 1990 by Hilton and Johnson \cite{AJWHilton1990}, and was also touched on in the survey paper by Woodall \cite{WRWoodall2001} in 2001.

For an edge-list assignment $L$ of $G$, a proper {\it L-edge-colouring} is a mapping $\phi$ from $E(G)$ to $\bigcup \limits_{e\in E}L(e)$ such that $\phi(e) \in L(e)$ for all $e \in E(G)$, and if $e, f \in E(G)$  have a vertex in common, then $\phi(e)\ne \phi(f)$. For a list assignment $L$ to $E(G)$ and a colour $\sigma \in \bigcup \limits_{e\in E(G)}L(e)$, let $\alpha'(\sigma,L,G)$ be the maximum number of independent edges in the subgraph $H$ of $G$ induced by those edges in $E(G)$ such that $\sigma \in L(e)$. Thus $\alpha'(\sigma,L,G)$ is the size of the largest set of edges of $G$ having $\sigma$ in their lists and which have pairwise no vertex in common; in other words $\alpha'(\sigma,L,G)$ is the {\it matching number} of the subgraph of $G$ induced by the edges with $\sigma$ in their lists. In the case when each edge has the same list assignment, the notation $\alpha'(\sigma,L,G)$ is shortened to  $\alpha'(G)$. 

Suppose that $E(G)$ has a proper L-edge colouring. Then for each $\sigma \in \bigcup \limits_{e\in E} L(e)$ , the set $S_{\sigma}$ of edges coloured  $\sigma$ is independent. Therefore, since each vertex is coloured,
\begin{equation*}
|E(G)|=\sum_{\sigma} |S_\sigma| \le\sum_{\sigma} \alpha'(\sigma,L,G), 
\end{equation*}

\noindent where the sums are taken over all $\sigma \in \bigcup \limits_{e\in E} L(e)$. Since every subgraph $H$ of $G$ is also properly coloured, the same inequality holds with $E=E(G)$ replaced with $E(H)$ and $\alpha'(\sigma,L,G)$ replaced by $\alpha'(\sigma,L,H)$.

Given a simple graph $G$ and a list assignment $L$ for $E(G)$, we say that $G$ satisfies {\it Hall’s edge-condition} if and only if, for each subgraph $H$ of $G$,

\begin{equation*}
|E(H)|\le \sum_{\begin{subarray}{c}
\sigma \in (\cup L(e): e \in E(H))
\end{subarray}} \alpha'(\sigma,L,H)\hspace{20pt} (*)
\end{equation*}

\noindent Note that (*) holds for each induced subgraph $H$ of $G$ if and only of it holds for each subgraph of $G$.

We define the {\it Hall edge number} or {\it Hall index} $h'(G)$ to be the smallest positive integer $l$ such that there is a proper L-edge-colouring of $G$ whenever $G$ and $L$ satisfy Hall’s edge-condition and $|L(e)|\ge l$ for all $e \in E(G)$.

We define the {\it Hall edge condition number}, or the {\it Hall condition index $s'(G)$} of $G$ to be the smallest integer $l$ such that $G$ and $L$ satisfy Hall’s Condition (*) whenever $|L(e)|\ge l$ for all $e \in E(G)$. Let $s_{o}{'}(G)$ be the smallest integer $l$ such that the assignment of $\{1,...,l\}$ to every edge satisfies Hall’s Condition.

\noindent {\it Theorem 14.} (Hilton and Johnson \cite{AJWHilton1999} and Johnson \cite{PDJohnsonJr1999})
\begin{equation*}
\begin{split}
s'(G)=s_{o}{'} & =max\Bigg\{{\Biggl\lceil\frac{|E(H)|}{\alpha'(H)}\Biggr\rceil}:\text{\it{H is an induced subgraph of G}} \Bigg\}\\
& =max\Bigg\{{\Biggl\lceil\frac{|E(H)|}{\alpha'(H)}\Biggr\rceil}:\text{\it{H is an induced subgraph of G with at least one edge}} \Bigg\}.
\end{split}
\end{equation*}

Some further facts about $c'(G),s'(G),\chi'(G)$ and $h'(G)$ are collected together in the following Theorem. (See Hilton and Johnson \cite{AJWHilton1999}, Vizing \cite{VGVizing1964} and Theorem 5).
\newpage
\noindent {\it Theorem 15. Let $G$ be a simple graph. Then}
\begin{enumerate}[label=\Alph*)]
\item { Either}
\begin{enumerate}[label=\arabic*)]
\item $\Delta(G)+1=c'(G)=\chi'(G)=s'(G)>h'(G)$, or
\item $\Delta(G)+1=c'(G)=\chi'(G)=h'(G)\ge s'(G)$, or
\item $\Delta(G)+2\ge c'(G)=h'(G)>\chi'(G)\ge s'(G)\ge \Delta(G)$\\
(so $c'(G)=max\{h'(G),s'(G)\}$ and $c'(G)=max\{h'(G),\chi'(G)\}$.
\end{enumerate}
\item $\Delta(G)+1\ge s'(G)\ge\Delta(G)$,\\
$\Delta(G)+2\ge h'(G)$.
\item $s'(G)-s'(G-e)\le 1$ if $\Delta(G-e)=\Delta(G)$,\\
$h'(G)-h'(G-e)$ can be arbitrarily large.
\item If $J$ is a subgraph of $G$ then $h'(J)\le h'(G)$ and $s' (J)\le s'(G).$
\end{enumerate}

The most interesting question about the choice index $c'(G)$ is whether the choice index conjecture is true, i.e. whether $c'(G)=\chi'(G)$ for all simple graphs $G$, or indeed, for all multigraphs. If it is not true, then (A 3) above has to be true for some graph $G$, and so either $\Delta(G)+2=c'(G)=h'(G),$ $ \Delta(G)+1 \ge \chi(G) \ge s'(G)=\Delta(G)$ or  $\Delta(G)+1=c'(G)=h'(G)$, $\chi'(G)=s'(G)=\Delta(G)$.  Unfortunately as a test whether some particular graph satisfies the choice index conjecture, this is not very easy to carry out, as it is usually very difficult to evaluate $h'(G)$.This aspect was looked at by Cropper and Hilton \cite{MMCropper2002}. They were unable to determine the choice index or the Hall index of $K_{2n}$.

We remark that the Overfull Conjecture of Chetwynd and Hilton \cite{Chetwynd1986Critical} can be expressed in the form:
	
\noindent {\it Conjecture 16}. Let $\Delta(G)>\frac{1}{3} |V(G)|$. Then $G$ is class 2 if and only if $s'(G)=\Delta(G)+1$.\\
(See Hilton and Johnson \cite{AJWHilton1999} and Dugdale, Eslahchi and Hilton \cite{JKDugdale2000})\\
In fact, in the context of simple graphs $G$, the parameter $s'(G)$ might be called the {\it overfullness} of $G$.

\section{The total Hall number and the total Hall Condition number of a simple graph.}

If $H$ is a subgraph of $G$ induced by $E(H)\cup V(H)$ then let $\alpha_T(H)$ be the largest independent set of vertices and edge of $G$, or, in other words, let $\alpha_T(H)$ be the largest set of vertices and edge of $G$ such that no two vertices are adjacent, no two edges are incident with the same vertex, and it contains no edge-vertex pair with the edge incident with the vertex. $\alpha_T(H)$ is called the {\it total independence number} of $H$.

Suppose that $E(G)\cup V(G)$ has a proper L-total colouring. For a total list assignment $L$ to $E(G)\cup V(G)$ and a colour $\sigma \in (\bigcup \limits_{e\in E} L(e)) \cup ( \bigcup \limits_{v\in V} L(v)$). Let $\alpha_T (\sigma,L,G)$ be the maximum number of independent edges and vertices in the subgraph $H$ of $G$ which contain $\sigma$ in their lists. Then, for each $\sigma \in (\bigcup \limits_{e\in E} L(e))\cup \bigcup \limits_{v\in V} L(v))$, the set $T_{\sigma}$ of edges and vertices is independent. Therefore, since each vertex and edge is coloured,

\begin{equation*}
|E(G)|\cup |V(G)| = \sum_{\begin{subarray}{c}
\sigma
\end{subarray}} |T_{\sigma}|
=\sum_{\begin{subarray}{c}
\sigma
\end{subarray}} \alpha_{T}{(\sigma,L,G),}
\end{equation*}

\noindent where the sum is taken over all $\sigma \in (\bigcup \limits_{e\in E} L(e)) \cup (\bigcup \limits_{v\in V} L(v))$.

\noindent Since every subgraph $H$ of $G$ is also properly coloured, the same inequality holds with $E(G)$ and $V(G)$ replaced by $E(H)$ and $V(H)$, and $\alpha_T(\sigma,L,G)$ by $\alpha_T (\sigma,L,H)$.

Given a simple graph $G$ and a total list assignment $L$ for $E(G)\cup V(G)$, we say that $G$ satisfies {\it Hall’s total condition} if and only if, for each induced subgraph $H$ of $G$,

\begin{equation*}
|V(H)|+|E(H)|\le \sum_{\begin{subarray}{c}
\sigma \in L(e)\cup L(v)\\
e \in E(H),v \in V(H)
\end{subarray}} \alpha_{T}(\sigma,L,H) \hspace{1cm} (**)
\end{equation*}

We define the {\it total Hall number},$h_T(G)$, of $G$ to be the smallest positive integer $l$ such that there is a proper L-total-colouring of $G$ whenever $G$ and $L$ satisfy Hall’s total condition (**) whenever $|L(e)\cup L(V)|\ge l$ for all $e \in E(G)$ and $v \in V(G)$.

We define the {\it total Hall condition number} $s_T(G)$ to be the smallest integer $l$ such that $G$ and $L$ satisfy Hall’s total condition (**) whenever $|L(e)|\ge l$ and $|L(v)|\ge l$ for all $e \in E(G)$ and  $v \in V(G)$. Let $s_{To} (G)$ be the smallest integer $l$ such that the assignment of $\{1,...,l\}$ to every edge and every vertex satisfies Hall’s total condition (**).

\noindent {\it Theorem 17 (Hilton and Johnson \cite{AJWHilton1999} and Johnson \cite{PDJohnsonJr1999}).}

\noindent {\it For a simple graph $G$}, 
\begin{equation*}
s_T(G)=s_{To}(G)=max\Bigg\{{\Biggl\lceil\frac{|V(H)+|E(H)|}{\alpha_T(H)}\Biggr\rceil}:\text{\it{H is an induced subgraph of G}} \Bigg\}
\end{equation*}

It is not known if $h_T(G)\ge h_T(G-e)$ for each edge $e \in E(G)$, but the following is true. If $H$ is a subgraph of $G$ induced by some subset of $V(G)\cup E(G)$, then $h_T(H)\le h_T(G)$. It is also not known if $h_T(G)-h_T (G-e)$ can be greater than $1$. It is true that $h_T (G)\ge max\{h(G),h'(G)\}$, where $h(G)$ is the Hall number of $G$ (defined in \cite{AJWHilton1999} and elsewhere).

Some further facts about $c_T(G),s_T(G),\chi_T(G)$ and $h_T(G)$ are collected together in the following theorem. (see Hilton and Johnson \cite{AJWHilton1999} and Theorem 11).

\noindent {\it Theorem 18. Let $G$ be a simple graph. Then,}
\begin{enumerate}[label=\Alph*)]
\item {\it Either}
\begin{enumerate}[label=\arabic*)]
\item $\Delta(G)+4\ge c_T(G)=\chi_T(G)=s_T(G)>h_T(G)$, {\it or}
\item $\Delta(G)+4\ge c_T(G)=\chi_T(G)=h_T(G)\ge s_T(G)\ge \Delta(G)+1$, {\it or}
\item $\Delta(G)+4\ge c_T(G)=h_T(G)>\chi_T(G)\ge s_T(G)\ge \Delta(G)+1$.
\end{enumerate}
\item $\Delta(G)+4\ge s_T(G)\ge \Delta(G)+1$,\\
$\Delta(G)+4\ge h_T(G)$,
\item {\it If $h$ is a subgraph of $G$ induced by some subset of $V(G)\cup E(G)$, then $h_T(H)\le h_T(G)$ and $s_T(H)\le s_T(G)$.}
\end{enumerate}
\newpage
\bibliographystyle{unsrt}  


\end{document}